\documentclass[11pt]{amsart}
\usepackage{fullpage}
\usepackage{amsmath, amssymb, amsthm}

\usepackage{graphicx}
\usepackage{color}
\usepackage[colorlinks=true,linkcolor=red,citecolor=blue,urlcolor=cyan]{hyperref}
\usepackage{algorithm}
\usepackage{algpseudocode}

\newtheorem{thm}{Theorem}[section]

\graphicspath{{./figs/}} %note trailing /

\newcommand{\ba}{\mathbf a}

\newcommand{\be}{\mathbf e}

\newcommand{\norm}[1]{\left\Vert#1\right\Vert}

\newcommand{\B}{\mathcal B}
\newcommand{\D}{\mathcal D}

\newcommand{\tF}{\tilde{F}}

\newcommand{\Lo}{\mathcal{L}}

\newcommand{\bbZ}{\mathbb Z}

\newcommand{\bzeta}{\boldsymbol{\zeta}}

\newcommand{\bkappa}{\boldsymbol\kappa}

\newcommand{\sgn}{\mathop{\mathrm{sgn}}}

\newcommand{\bx}{\mathbf x}

\usepackage[backend=biber,giveninits=true,sorting=nyt,style=alphabetic,natbib=true,maxcitenames=2,maxbibnames=10,url=false,doi=true,backref=false]{biblatex}
\renewbibmacro{in:}{\ifentrytype{article}{}{\printtext{\bibstring{in}\intitlepunct}}}

\begin{document}

\title{A Semi-definite Optimization Method for Maximizing the Shared Band Gap of Topological Photonic Crystals}

\author{Chiu-Yen Kao}
\address{Department of Mathematical Sciences, Claremont McKenna College, Claremont, CA }
\email{ckao@cmc.edu}
\thanks{Chiu-Yen Kao acknowledges partial support from NSF DMS 2208373.}

\author{Junshan Lin}
\address{Department of Mathematics and
Statistics, Auburn University, Auburn, AL 36849 }
\email{jzl0097@ auburn.edu}
\thanks{Junshan Lin acknowledges partial
supported by the NSF grants DMS-2011148 and DMS-2410645.}

\author{Braxton Osting}
\address{Department of Mathematics, University of Utah, Salt Lake City, UT}
\email{osting@math.utah.edu}
\thanks{Braxton Osting acknowledges partial support from NSF DMS-1752202.}

% 35P05 General topics in linear spectral theory for PDEs
% 35B10 Periodic solutions to PDEs
% 35Q93 PDEs in connection with control and optimization
% 49R05 Variational methods for eigenvalues of operators
%(should also be assigned at least one other classification number in Section 49) [See also 47A75]
%65K10 Numerical optimization and variational techniques
% 65N25 Numerical methods for eigenvalue problems for
%boundary value problems involving PDEs


\date{\today}

\begin{abstract}
Topological photonic crystals (PCs) can support robust edge modes to transport electromagnetic energy in an efficient manner. Such edge modes are the eigenmodes of the PDE operator for a joint optical structure formed by connecting together two photonic crystals with distinct topological invariants, and the corresponding eigenfrequencies are located in the shared band gap of two individual 
photonic crystals. 
This work is concerned with maximizing the shared band gap of two photonic crystals with different topological features in order to increase the bandwidth of the edge modes. We develop a semi-definite optimization framework for the underlying optimal design problem, which enables efficient update of dielectric functions at each time step while respecting symmetry constraints and, when necessary, the constraints on topological invariants. At each iteration, we perform sensitivity analysis 
of the band gap function and the topological invariant constraint function to linearize the optimization problem and solve a convex semi-definite programming (SDP) problem efficiently.
Numerical examples show that the proposed algorithm is superior in generating optimized optical structures with robust edge modes.

% Designing robust surface (edge) modes are of great importance in photonic crystals due to their emerging applications in waveguides, lasers, and information transmission. In this work, we propose computational approaches to create two photonic crystals that have large common bandgaps and are of distinct Berry curvatures. These two photonic crystals can then be glued together to generate edge modes propagating along the interface in between. The sensitivity analysis of bandgaps and Berry curvatures are carried out. A semi-definite optimization framework enables efficient updating of dielectric distributions at each time step while respecting symmetry constraints and, when necessary, Berry curvatures constraints. We show that our algorithms are superior in generating optimized structures with robust edge modes.    
\end{abstract}

\maketitle

\section{Introduction}

Photonic crystals with band gaps can be used to confine optical waves and have important applications in transporting wave energy \cite{Joannopoulos_2011}.
The recent development in topological insulators (cf. \cite{hasan2010colloquium, bernevig2013topological, qi2011topological}) opens up new avenues for electromagnetic wave confinement in photonic crystals, when it was realized in \cite{raghu2008analogs} that topological band structures are a ubiquitous property of waves in periodic media, regardless of the classical or quantum nature of the waves. In particular, the concepts in topological insulators can be extended in parallel to acoustic, electromagnetic and mechanical waves in periodic media, and topological wave insulators can be built to manipulate the classical wave in the same way as solids modulating electrons \cite{khanikaev2013photonic, lu2014topological, ozawa2019topological, yang2015topological}. 

Wave confinement in topological insulators is achieved by gluing together two periodic media with distinct topological invariants. The topological phase transition at the interface of two media gives rise to the edge modes, which propagate parallel to the interface but localize in the direction transverse to the interface. Such edge modes are topologically protected and are stable against fabrication imperfections, which make them attractive for wave energy transportation in realistic applications. We refer to \cite{khanikaev2013photonic, lu2014topological, ozawa2019topological} and references therein for the experimental investigation of edge modes in photonic crystals, and \cite{ammari2023edge, ammari2020topologically, bal2019continuous, drouot2020edge, fefferman2017topologically, lee2019elliptic, lin2021mathematical, qiu2023mathematical} and the references therein for the mathematical studies of edge modes for several elliptic operators.

The topologically protected edge mode, which is responsible for the transportation of wave energy along the material interface, attains an eigenfrequency in the spectral band gap
of two bulk media glued together. Therefore, it is desirable to increase the width of the band gap to allow for the existence of edge modes within a larger bandwidth and the operation of energy transportation in a more stable manner.
For photonic crystals, there has been extensive work on maximizing the length of spectral band gaps.
Among one-dimensional photonic crystals, the Bragg structure can be proven to have the largest first spectral gap-to-midgap ratio \cite{osting2012bragg}. 
The two-dimensional photonic crystal that maximizes the $m$-th spectral gap-to-midgap ratio is unknown analytically, but computational methods have been developed for its approximation. 
In particular, gradient-based methods were developed in 
\cite{dobson1999maximizing,cox2000band,kao2005maximizing,sigmund2008geometric}. Methods based on semi-definite programming were developed in 
\cite{Men_2010,Men_2014,Kao_2019} 
and shown to better handle multiplicities in the eigenvalues, improving convergence. 
For topological photonic crystals, the available computational strategies of manually adding/removing materials from the existing structures to open the gap near Dirac points or choosing the parameter that yields the largest gap value among a set of available parameters are inefficient \cite{nanthakumar2019inverse, nussbaum2021inverse}. The topology optimization developed in \cite{chen2019inverse,chen2021creating,chen2022inverse,chen2022inversely} relies on choosing the band gap interval first and  then exciting the desired modes with suitable current sources at the band gap boundary to achieve desired topological phases, but there is no guarantee that such a prechosen band gap would be the maximum. 

In this work, we develop an optimization framework to maximize the shared band-gap width of two photonic crystals directly,  subject to the PDE constraints in the form of the eigenvalue problem 
%\eqref{e:twistHelm} 
and the desired topological phase constraints. 
To solve the optimization problem, we extend the semi-definite programming  approach, while respecting symmetry constraints. 
Furthermore, we perform the sensitivity analysis of the topological phase with respect to the medium coefficient of the photonic crystal; we linearize the Berry curvature and impose the topological phase constraint effectively at each iteration. Numerical examples show that the algorithm converges fast when applied to the optimization problem, and the optimized photonic structures attain much larger band gaps compared to the initial structures; see Figure \ref{Fig:square_lattice_opt_PCs} for results on a square lattice and Figures \ref{Fig:tri_C3_opt_PCs} and  \ref{Fig:tri_C6_opt_PCs} for results on hexagonal lattices with C3 and C6 symmetry, respectively. The obtained photonic crystals possess much larger shared band gap than the initial structures
%and their joint structures support the propagation of edge modes.  
We further demonstrate that the joint structure, formed by gluing two optimized photonic crystals with desired topological indexes, indeed support robust transport of edge modes over a wide frequency range. 

The rest of the paper is organized as follows. 
In Section~\ref{s:MathModel}, we formulate the spectral problems and introduce the underlying topological invariant and topologically protected edge mode in photonic crystals. 
In Section~\ref{s:CompMethods}, we formulate the optimization problem and present the semi-definite programming method to maximize the shared band gap of two  photonic crystals with distinct topological invariants. 
%To enforce the constraints on the topological invariant effectively, we also perform the sensitivity analysis of the Berry curvature with respect to the medium coefficient of the photonic crystal. 
Numerical results for both square and hexagonal lattices are demonstrated in Section~\ref{s:CompResults}. We conclude our paper and discuss some potential future work in Section~\ref{s:Disc}. 

%\clearpage 
\section{Mathematical model for topological photonic crystals}
\label{s:MathModel}

\subsection{Time-harmonic waves in photonic crystals}
We consider materials with spatially periodic electric permittivity $\varepsilon(\bx)$ and magnetic permeability $\mu(\bx)$,  which are called \emph{photonic crystals} \cite{Joannopoulos_2011}. 
The propagation of time-harmonic electromagnetic waves in photonic crystals is governed by the time-harmonic Maxwell's equations, 
\begin{subequations}
\label{e:Maxwell}
\begin{align}
& \nabla \times \mathbf E = \imath \omega \mu \mathbf H, \\
& \nabla \times \mathbf H = - \imath \omega \varepsilon \mathbf E,
\end{align}
\end{subequations}
where $\mathbf E$ and $\mathbf H$ are the electric and magnetic fields respectively and  $\omega$ is the angular frequency. 
Let us write 
$\varepsilon(\bx) = \varepsilon_0 \varepsilon_r(\bx)$, 
$\mu(\bx) = \mu_0 \mu_r(\bx)$, 
where $\varepsilon_r$, $\varepsilon_0$, $\mu_r$,  and $\mu_0$ are the relative permittivity, vacuum permittivity, relative permeability, and vacuum permeability, respectively.
Here we consider the non-magnetic material by assuming that $\mu_r \equiv 1$.
Due to the scaling invariance of the Maxwell's equations, without loss of generality, we set the wave speed $c := \left( \varepsilon_0 \mu_0 \right)^{-1/2}=1$. 
Furthermore, we consider the two-dimensional problem by assuming that the material is invariant along the $z$-direction and the electromagnetic wave attains the transverse electric (TE) polarization with the electric field $\mathbf E = (0,0,E_z)$.  Writing $\varepsilon_r = \varepsilon$,  then $E_z$ satisfies the following Helmholtz equation
\begin{equation}
\label{e:Helm}
 - \Delta E_z(\bx) = \lambda \, \varepsilon(\bx) E_z(\bx), \qquad \qquad \bx \in \mathbb R^2,
\end{equation}
where $\lambda = (\omega/c)^2 = \omega^2$.

For a periodic lattice $\Lambda$ in $\mathbb{R}^2$
given by $ \Lambda := \bbZ \be_1 \oplus  \bbZ \be_2$, wherein $\be_1$ and $\be_2$ are the lattice vectors,
we consider $\varepsilon$ to be $\Lambda$-periodic, i.e., $\varepsilon(\bx+\be) = \varepsilon(\bx)$ for $\bx \in \mathbb R^2$ and $\be \in \Lambda$. 
The Floquet-Bloch transform (cf. \cite{Kuchment_1993}) can be applied to investigate the spectrum of the PDE operator in \eqref{e:Helm}.  To this end, let $\B \subset \mathbb R^2$ be the \emph{Brillouin zone}, which can be taken as the Voronoi cell of the origin in the reciprocal lattice $\Lambda^*$. 
For each $\bkappa \in \B$, we seek the Floquet-Bloch mode for \eqref{e:Helm} in the form of $E_z(\bx,\bkappa) = e^{\imath \bkappa\cdot \bx} \phi(\bx,\bkappa)$, where $\phi(\cdot,\bkappa)$ is a $\Lambda$-periodic function satisfying 
\begin{subequations}
\label{e:twistHelm}
\begin{align}
\label{e:twistHelmA}
& - (\nabla + \imath \bkappa) \cdot (\nabla + \imath \bkappa)  \phi(\bx) = \lambda \,\varepsilon(\bx) \phi(\bx)  && \bx \in \mathbb R^2  \\ 
\label{e:twistHelmB}
& \phi(\bx + \be) = \phi(\bx) && \be \in \Lambda.
\end{align}
\end{subequations}
By periodicity, we can restrict $\bx \mapsto \phi(\bx,\bkappa)$ to the fundamental periodic cell  $\D:=\{ \ell_1 \be_1 + \ell_2 \be_2; 0 \le \ell_1 \le 1, 0 \le \ell_2 \le 1\}$.
Let $L_\varepsilon^2(\D)$ be the Hilbert space of square integrable functions equipped with the inner product
\[
\langle \phi,\psi\rangle =\int_\D \varepsilon(\bx)\phi(\bx) \overline{\psi}(\bx) \, d\bx.
\]
We also define the  Hilbert spaces of periodic functions 
$$L_{\varepsilon,p }^2(\D):=\{ \phi(\bx)\in L_\varepsilon^2(\D)\colon  \, \, \phi(\bx+\be) = \phi(\bx), \, \forall \be \in \Lambda \}$$ and
$$
H_{\varepsilon,p}^{1}(\D):=\{ \phi(\bx)\in L_{\varepsilon,p}^2(\D)\colon \, \, \partial_i\phi(\bx)\in L_{\varepsilon,p}^2(\D), \ i=1,2 \}. $$
For each $\bkappa\in \B$, the differential operator $ A: = - \frac{1}{\varepsilon} (\nabla + \imath \bkappa) \cdot  (\nabla + \imath \bkappa) $ restricted to $L_{\varepsilon,p}^2(\D)$ is self-adjoint and it has a discrete positive spectrum with eigenvalues ordered by (cf. \cite{Kuchment_1993}) 
$$
\lambda_1(\bkappa) \leq \lambda_2(\bkappa) \leq \cdots \to \infty .
$$
The \emph{dispersion relation} (Bloch variety) consists of $(\bkappa,\lambda) \in \B \times \mathbb R$ such that the eigenvalue problem \eqref{e:twistHelm} attains a nontrivial solution. 
It is convenient to decompose the dispersion relation into spectral bands, 
$\lambda_m(\bkappa)$, which are $\Lambda^*$-periodic and can be considered over the  Brillouin zone $\B$. 

%\subsection{Band-gap for photonic crystal}
We say that that there is a \emph{spectral gap between the $m$-th and $(m+1)$-th spectral bands} if 
$\min_{\bkappa \in \B} \lambda_{m+1}(\bkappa) \ > \ 
\max_{\bkappa \in \B} \lambda_m(\bkappa)$. 
When there is a spectral gap between the $m$-th and $(m+1)$-th spectral bands, we will be interested in the length of the spectral gap, 
$$
\min_{\bkappa \in \B} \ \lambda_{m+1}(\bkappa)
\ - \ 
\max_{\bkappa \in \B} \ \lambda_m(\bkappa). 
$$
In order to better characterize the length of the gap so that it is independent of scaling, it is common to consider the
\emph{gap-to-midgap ratio} 
\begin{equation}
\label{e:Gi}
G_m(\varepsilon) := 2 \frac{ \min_{\bkappa \in \B} \omega_{m+1}(\varepsilon) - \max_{\bkappa \in \B} \omega_{m}(\varepsilon)}{\min_{\bkappa \in \B} \omega_{m+1}(\varepsilon) + \max_{\bkappa \in \B} \omega_{m}(\varepsilon)}
\end{equation} 
where $\omega_{m+1}^2=\lambda_{m+1}$
and $\omega_{m}^2=\lambda_{m}$.

\subsection{Topological invariant of photonic crystals}
The topological phase of photonic crystal is defined via the periodic part of the Bloch mode that satisfies \eqref{e:twistHelm}  and  $\norm{\phi(\bkappa; \cdot)}_{L_\varepsilon^2(\D)}=1$. Here and henceforth, we assume the $k$-th band $\lambda_k(\bkappa)$ of the photonic crystal is a simple eigenvalue for each $\bkappa\in\B$. Then the \textit{Berry phase} associated with the $k$-th band is given by the following integral along the boundary of the Brillouin zone \cite{bernevig2013topological, vanderbilt2018berry}:
\begin{equation}\label{eq:Berry_phase1}
\theta_k := \oint_{\partial \B} \bzeta_k(\bkappa) \cdot d\bkappa, 
\end{equation}
where
$\bzeta_k(\bkappa) := i \big\langle  \phi_k(\cdot, \bkappa), \nabla_{\bkappa} \phi_k(\cdot, \bkappa) \big\rangle _{L^2(\mathbb R^2 / \Lambda)}$ is the \emph{Berry connection}.  
%In the above, $\langle  \cdot, \cdot \rangle $ denotes the $L^2-$inner product over the fundamental cell of the periodic media. 
Alternatively, using the Stoke's theorem, the Berry phase is written as
\begin{equation}\label{eq:Berry_phase2}
\theta_k = \int_\B  F_k(\bkappa) \, d\bkappa,
\end{equation}
where
\begin{equation}\label{eq:Berry_curvature}
   F_k(\bkappa) := \partial_{\kappa_x}  \bzeta_{k}^{(2)}(\bkappa) -  \partial_{\kappa_y}  \bzeta_{k}^{(1)}(\bkappa) 
\end{equation}
is the \emph{Berry curvature}, and $\bzeta_{k}^{(j)}$ denotes the $j$-th component of the vector $\bzeta_{k}$.
It can be shown that the number
\begin{equation}\label{eq:Chern_number}
C_k :=  \frac{\theta_k}{2\pi} = \frac{1}{2\pi} \int_\B F_k(\bkappa) \, d\bkappa 
\end{equation}
is an integer, which is called the \textit{Chern number} and is a \textit{topological invariant} of the medium \cite{bernevig2013topological, vanderbilt2018berry}.
Assume that there is a gap between the $m$-th and $(m+1)$-th band of the spectrum, the gap Chern number is defined by $C_{g,m} = \sum_{j=1}^m C_j$.

When the Chern number is trivial,
other types of topological invariants can be defined. 
For instance, one can use \textit{valley Chern number} as a topological invariant if the system attains a valley degree of freedom in the Brillouin zone.
This applies to periodic media with broken inversion symmetry when the Berry curvature at two symmetric points  $\bkappa_1$ and $\bkappa_2$ in the Brillouin zone attain opposite signs.
By integrating the Berry curvature over the two regions $\B(\bkappa_1)$ and $\B(\bkappa_2)$ near the symmetric valley points,
the valley Chern number $C_k^{(\bkappa_1, \bkappa_2)}$ for the $k$-th band can be defined as 
\begin{equation}\label{eq:valley_Chern_number}
C_k^{(\bkappa_1, \bkappa_2)}:= \sgn \left(C_k^{(\bkappa_1)} -  C_k^{(\bkappa_2)} \right), \qquad \mbox{where} \ 
C_k^{(\bkappa_j)} := \int_{\B(\bkappa_j)} F_k(\bkappa) \, d\bkappa, 
\  j \in [2]\footnote{Here, we use the notation $[n] := \{1,\ldots,n\}$ for $n \in \mathbb N \setminus\{0\}$.}.
\end{equation}
The valley Chern number $C_k^{(\bkappa_1, \bkappa_2)}$ above is a topological invariant and the corresponding photonic crystal mimics the so-called quantum valley Hall effect studied in electron models \cite{dong2017valley, kim2021multiband, Ma_2016, wong2020gapless}. Similar to the gap Chern number, the gap valley Chern number between the $m$-th and $(m+1)$-th band of the spectrum is defined by
\begin{equation}\label{eq:gap_valley_Chern_number}
 C_{g,m}^{(\bkappa_1, \bkappa_2)}:= \sgn \left(C_{g,m}^{(\bkappa_1)} -  C_{g,m}^{(\bkappa_2)} \right),   
\end{equation}
% $C_{g,m}^{(\bkappa_1, \bkappa_2)}:= \sgn \left(C_{g,m}^{(\bkappa_1)} -  C_{g,m}^{(\bkappa_2)} \right)$, 
where $C_{g,m}^{(\bkappa_j)} := \int_{\B(\bkappa_j)} F_{g,m}(\bkappa) \, d\bkappa$, and the gap Berry curvature $F_{g,m}(\bkappa)=\sum_{j=1}^m F_j(\bkappa)$.

Another well-known topological invariant is the so-called spin Chern number $C_k^{s}$, which is introduced for the spin Hall effect of electron models \cite{wu2015scheme}. The spin Chern number $C_k^{s}$ can be computed by separating the two distinct spin eigenmodes and performing the Chern number calculations.
% for each using the formula \eqref{eq:Chern_number}, and similar to the gap Chern number, the gap spin Chern number is computed by summing up $C_k^{s}$ for the first $m$ bands \cite{khanikaev2013photonic}.
There also exist other topological phases defined in different settings, including Zak phase, Wilson loop, etc \cite{alexandradinata2014wilson, liu2018topological}.

\subsection{Topologically protected edge mode}
 When two photonic crystals with the permittivity function $\varepsilon(\bx)=\varepsilon_1(\bx)$ and $\varepsilon(\bx)=\varepsilon_2(\bx)$ attain distinct topological invariants and are glued together, a topological phase transition takes place at their interface. Such phase transition gives rise to \textit{edge modes} propagating parallel to the interface but localized in the direction transverse to the interface.
Mathematically, for each parallel wave vector $\bkappa_\parallel$, 
the edge mode $\psi(\bkappa_\parallel; \bx)$ solves the following eigenvalue problem for $\lambda$ located in the shared spectral band gap of the two bulk periodic media:
\begin{align}
& -\Delta \psi(\bkappa_\parallel; \bx) = \lambda \,\varepsilon(\bx) \psi(\bkappa_\parallel; \bx) \quad \mbox{for} \; \bx\in\mathbb{R}^2, \label{eq:edge_mode_PDE} \\
& \psi(\bkappa_\parallel; \bx+\ba)=e^{\imath \bkappa_\parallel \cdot \ba}\psi(\bkappa_\parallel; \bx) \hspace*{0.6cm}( \psi \; \mbox{propagates along} \; \ba ), \label{eq:edge_mode_qp_bnd} \\
& \psi(\bkappa_\parallel; \bx) \to 0  \;  \mbox{as} \;  |\bx \cdot \ba^\perp| \to\infty \;   \hspace*{0.5cm} ( \psi  \; \mbox{localized along}  \; \pm\ba^\perp). \label{eq:edge_mode_decay}
\end{align}
In the above, $\ba$ is the lattice vector along the interface direction, and $\ba^\perp$ is the  unit vector orthogonal to $\ba$. In addition, the permittivity function $\varepsilon(\bx)=\varepsilon_1(\bx)$ and $\varepsilon(\bx)=\varepsilon_2(\bx)$ for $\bx\cdot\ba_\perp>0$ and  $\bx\cdot\ba_\perp<0$ respectively. The spectra of the edge modes are topologically protected, meaning their existence is guaranteed
by the topological nature of two bulk materials, which allows for fabrication imperfections in the bulk media. In general, the number of edge modes is equal to the difference of the bulk topological invariants across the interface, which is known as the \textit{bulk-edge correspondence} \cite{asboth2016short, bernevig2013topological}.

\section{Computational algorithm for the band-gap optimization} \label{s:CompMethods}

\subsection{Formulation of the optimization problem}\label{ss:forOptim}
Here, we will consider the problem of designing two periodic media with distinct topological invariants so that they may \emph{optimally} support topologically protected edge modes. Since the edge mode must lie in the shared spectral gap, we formulate the problem by maximizing the shared gap-to-midgap ratio. 
%Denote the collection of medium coefficients by 
%$\zeta = (\varepsilon_1,\varepsilon_2)$. 

We first define the \emph{admissible class of pairs of permittivities}, denoted by $\mathcal{A}$. 
Fix 
$\varepsilon_+ > \varepsilon_- > 0$, 
$m \in \mathbb N$,  
$n_i \in \{+1,-1\}$, $i\in [2]$. 
We say $(\varepsilon_1,\varepsilon_2)$ is an admissible pair of permittivities and write $(\varepsilon_1,\varepsilon_2) \in \mathcal{A}$ if it satisfies the following three conditions: 
\begin{enumerate}
\item[(A)] $\varepsilon_i\in L^\infty(\mathbb R^2)$, $i\in[2]$  is piecewise continuous and $\Lambda$-periodic,
\item[(B)] $\varepsilon_i$, $i \in [2]$  satisfies the pointwise bounds 
$$
\varepsilon_- \leq \varepsilon_i(\bx) \leq \varepsilon_+, 
\qquad \bx \in \D, \ i \in [2],  
$$
\item[(C)] 
The topological indices for the $m$-th band are fixed,
$$
%\label{e:TopoConstraint}
\mbox{Idx}_m(\varepsilon_i) = n_i, \qquad \qquad i \in [2]. 
$$ 
% the $m$-th valley Chern number for each material is fixed, 
% $$
%\label{e:TopoConstraint}
% C_m^{(\bkappa_1, \bkappa_2)}(\varepsilon_i) = n_i, \qquad \qquad i \in [2]. 
% $$
\end{enumerate}
Here $\mbox{Idx}_m$ represents the underling topological indices of the photonic crystal as discussed in Section 2. The Chern number is a natural topological index when the time-reversal symmetry of the photonic crystal is broken \cite{Joannopoulos_2011}. For the problems considered in this work, where the time-reversal symmetry is preserved and we optimize the electric permittivity function $\varepsilon_i(\bx)$ in the spatial domain for photonic crystals that mimics either the quantum valley or spin Hall effect, the topological invariant used is either gap valley Chern number or the gap spin Chern number.

We consider the smallest eigenvalue for the $(m+1)$-th band among the two materials, $\lambda_u= \lambda_u(\varepsilon_1,\varepsilon_2)$, and 
the largest eigenvalue for the $m$-th band among the two materials, $\lambda_l = \lambda_l(\varepsilon_1,\varepsilon_2)$:
\begin{align*}
\lambda_u(\varepsilon_1,\varepsilon_2) &:= \min_{i \in [2]} \ 
\min_{\bkappa \in \B} \ \lambda_{m+1} (\bkappa; \varepsilon_i), \\
\lambda_l(\varepsilon_1,\varepsilon_2) &:= \max_{i \in [2]} \ 
\max_{\bkappa \in \B} \ \lambda_{m} (\bkappa; \varepsilon_i). 
\end{align*}
The \emph{shared gap-to-midgap ratio} for the spectral gap between the $m$-th and $(m+1)$-th spectral bands is then defined by 
\begin{equation}
\label{e:J}
J_m(\varepsilon_1,\varepsilon_2) := 2 \frac{ \lambda_u(\varepsilon_1,\varepsilon_2) - \lambda_l(\varepsilon_1,\varepsilon_2)}{\lambda_u(\varepsilon_1,\varepsilon_2) + \lambda_l(\varepsilon_1,\varepsilon_2)}. 
\end{equation}
With $m$ fixed, we then consider the PDE-constrained optimization problem of maximizing the shared gap-to-midgap ratio among the admissible class of material coefficients 
\begin{align}
\label{e:OptProb}
\max_{(\varepsilon_1,\varepsilon_2)\in \mathcal{A}} \ &  J(\varepsilon_1,\varepsilon_2). 
\end{align}
In \eqref{e:OptProb}, it is understood that the eigenvalues $\lambda_{m} (\bkappa; \varepsilon_i)$ which define $J(\varepsilon_1,\varepsilon_2)$, satisfy \eqref{e:twistHelm} for the material coefficient $\varepsilon_i$.

In physical setting, the \emph{shared gap-to-midgap ratio in frequency} is used frequently due to the independence of scaling. We thus also define
\begin{equation}
\label{e:G}
G_m(\varepsilon_1,\varepsilon_2) := 2 \frac{ \omega_u(\varepsilon_1,\varepsilon_2) - \omega_l(\varepsilon_1,\varepsilon_2)}{\omega_u(\varepsilon_1,\varepsilon_2) + \omega_l(\varepsilon_1,\varepsilon_2)} 
\end{equation}
where $\omega_u^2=\lambda_u$ and $\omega_l^2=\lambda_l$. As $x \mapsto \max \{ 0, x^2\}$ is a monotone increasing function, $J_m$ and $G_m$ reach their maximums for the same pair of $(\varepsilon_1,\varepsilon_2)$. In the following, we will formulate a semi-definite optimization approach for maximizing  $J(\varepsilon_1,\varepsilon_2)$ to avoid using the square notation.

\subsection{Semi-definite programming optimization approach}
We reformulate  \eqref{e:OptProb} using a semi-definite programming optimization approach; see \cite{Men_2010,Kao_2019}. Introducing the parameters, $\alpha,\beta > 0$, \eqref{e:OptProb}  can be equivalently written as
\begin{subequations}
\label{e:OptProb2}
\begin{align}
\max_{\varepsilon_i, \alpha,\beta} \quad &  2 \frac{\beta - \alpha}{\alpha + \beta}  \\
\textrm{subject to} \quad 
& \beta \leq \lambda^{(i)}_{m+1}(\bkappa), && \bkappa \in \B, 
\ i \in [2] \\
& \alpha \geq \lambda^{(i)}_{m}(\bkappa), && \bkappa \in \B, \  i \in [2]\\
& - (\nabla + \imath \bkappa) \cdot (\nabla + \imath \bkappa) \  \phi_m^{(i)}(\bx) = \lambda^{(i)}_m \,\varepsilon_i(\bx) \phi_m^{(i)}(\bx)  && \bx \in \D, \ \bkappa \in \B, \ i\in [2]\\
& - (\nabla + \imath \bkappa) \cdot (\nabla + \imath \bkappa) \  \phi_{m+1}^{(i)}(\bx) = \lambda^{(i)}_{m+1} \,\varepsilon_i(\bx) \phi_{m+1}^{(i)}(\bx)  && \bx \in \D, \ \bkappa \in \B, \ i\in [2]\\
%& -\nabla \cdot \left( \rho_i^{-1}(x) \nabla \phi^{(i)}_j(x) \right)  = \lambda^{(i)}_j \eta_i(x) \phi^{(i)}_j(x), 
%&& x\in \mathbb R^2/\Lambda, \ j\in \mathbb N, \ i\in[2] 
& (\varepsilon_1,\varepsilon_2) \in \mathcal{A}. 
\end{align}
\end{subequations} 
The equivalence follows from the observations that for $\alpha,\beta>0$, 
$\alpha \mapsto \frac{\beta - \alpha}{\alpha + \beta}$ 
is decreasing  and 
$\beta \mapsto \frac{\beta - \alpha}{\alpha + \beta}$ 
is increasing. 
In the optimal case, we have $\alpha = \lambda^{(i)}_{m}(\bkappa)$ for some $ \bkappa \in \B$ and $i \in [2]$ as well as 
$\beta = \lambda^{(i)}_{m+1}(\bkappa)$ for some $ \bkappa \in \B$ and $i \in [2]$. 
To further reformulate \eqref{e:OptProb2}, we use the following SDP formulation of the min-max theorem. 
\begin{thm} \label{t:SDP-MM}
Let $\varepsilon \in L^\infty(\mathbb R^2)$ be piecewise continuous, positive, and $\Lambda$-periodic, $\bkappa \in \B$, and write $A = - \frac{1}{\varepsilon} (\nabla + \imath \bkappa) \cdot (\nabla + \imath \bkappa)$. For each $\bkappa\in\B$, $A$ is a self-adjoint operator on $L_{\varepsilon,p}^2(\D)$  and its eigenvalues are ordered by $\lambda_1\leq \lambda_2 \leq \cdots \to \infty $. Let $\mathcal{P}_k$ denote the collection of rank $k$ orthogonal projection operators and denote the identity operator by $\textrm{I}$.
Then 
\begin{subequations}
\label{e:MM1}
\begin{align}
\lambda_k  
\quad = \quad 
\inf_{\Pi \in \mathcal{P}_k}  \ & \lambda \\ 
\label{e:MM1b}
\textrm{s.t.} \ & \Pi \ (A - \lambda I ) \ \Pi \preceq 0. 
\end{align}
\end{subequations}
and 
\begin{subequations}
\label{e:MM2}
\begin{align}
\lambda_{k+1}  
\quad = \quad 
\sup_{\Pi \in \mathcal{P}_{k}}  \ & \lambda \\ 
\textrm{s.t.} \ & (I-\Pi) \ (A - \lambda I) \ (I - \Pi) \succeq 0. 
\end{align}
\end{subequations}
Furthermore, letting $\phi_j$ be normalized eigenfunctions of $A$ corresponding to $\lambda_j$, 
the extrema in \eqref{e:MM1} and \eqref{e:MM2} are each attained by the orthogonal project operator given by 
$\Pi \psi(\bx) = 
\sum_{j \in [k]} \left\langle  \phi_j (\cdot), \psi \right\rangle  \phi_j(\bx)$.
\end{thm}

% \begin{thm} \label{t:SDP-MM}
% Let $A \in \mathbb C^n$ be a Hermitian matrix with eigenvalues $\lambda_1\leq \lambda_2 \leq \cdots \leq \lambda_n$ and let $k \in [n]$. 
% Let $\mathcal{P}_k$ denote the set of rank $k$ orthogonal projection matrices and denote the identity operator by $\textrm{I}$.
% Then 
% \begin{subequations}
% %\label{e:MM1}
% \begin{align}
% \lambda_k  
% \quad = \quad 
% \min_{\Pi \in \mathcal{P}_k}  \ & \lambda \\ 
% %\label{e:MM1b}
% \textrm{s.t.} \ & \Pi \ (A - \lambda I ) \ \Pi \preceq 0. 
% \end{align}
% \end{subequations}
% and 
% \begin{subequations}
% %\label{e:MM2}
% \begin{align}
% \lambda_k  
% \quad = \quad 
% \max_{\Pi \in \mathcal{P}_{k}}  \ & \lambda \\ 
% \textrm{s.t.} \ & (I-\Pi) \ (A - \lambda I) \ (I - \Pi) \succeq 0. 
% \end{align}
% \end{subequations}
% Furthermore, letting $u_j$ be normalized eigenvectors of $A$ corresponding to $\lambda_j$, 
% the extrema in \eqref{e:MM1} and \eqref{e:MM2} are each attained by $\Pi = \sum_{j = 1}^k u_j \otimes u_j$. 
% %and the minimum in \eqref{e:MM2} is attained for $\Pi = I - \sum_{i = 1}^k u_i u_i^\ast = \sum_{i=k+1}^n u_i u_i^\ast$. 
% \end{thm}
\begin{proof}
The constraint \eqref{e:MM1b} is equivalent to 
$$
\frac{\int_{\D} \varepsilon u^*(\bx) A u(\bx) \ d\bx  }{\int_D \varepsilon u^*(\bx) u(\bx) \ d \bx} \leq \lambda,  
\qquad u \in \textrm{Ran} (\Pi) \setminus \{0\}.  
$$
Thus,  \eqref{e:MM1} can be written as 
$
\lambda_k  
=
\inf_{\Pi \in \mathcal{P}_k}  \ \ 
\sup_{u \in \textrm{Range} (\Pi)  \setminus \{0\} } \ \  \frac{\int_{\D} \varepsilon u^* A u \ d\bx  }{\int_D \varepsilon u^* u \ d \bx}$, which is the min-max variational characterization of $\lambda_k$ for a self-adjoint operator \cite[Ch.12]{Lieb_2001}.  
\eqref{e:MM2} is similarly proven using the max-min  characterization of $\lambda_{k+1}$. 
\end{proof}
Using Theorem~\ref{t:SDP-MM}, we can equivalently rewrite \eqref{e:OptProb2} as  

\begin{subequations}
\label{e:OptProb3}
\begin{align}
\max_{ \varepsilon_i, \alpha,\beta} \quad &  2 \frac{\beta - \alpha}{\alpha + \beta}  \\
\textrm{subject to} \quad 
& \Pi^{i,\bkappa} \left( - (\nabla + \imath \bkappa) \cdot (\nabla + \imath \bkappa)  - \alpha \varepsilon_i  \textrm{I} \right) \Pi^{i,\bkappa} \preceq  0, 
&& \bkappa \in \B_i, \  i \in [2] \\
& (\textrm{I} - \Pi^{i,\bkappa}) \left( - (\nabla + \imath \bkappa) \cdot (\nabla + \imath \bkappa)  - \beta \varepsilon_i  \textrm{I} \right) ( \textrm{I} - \Pi^{i,\bkappa})  \succeq  0, 
&& \bkappa \in \B_i, \  i \in [2] \\
& \textrm{$\Pi^{i,\bkappa}$ is a rank $m$ orthogonal projection}  
&& \bkappa \in \B_i, \  i \in [2] \\
& (\varepsilon_1,\varepsilon_2) \in \mathcal{A}. 
\end{align}
\end{subequations}
Introducing $\theta \geq 0$ and making the change of variables 
$\tilde \alpha = \theta \alpha^{-1}$, 
$\tilde \beta = \theta \beta^{-1}$, 
and
$\tilde \varepsilon_i = \theta \varepsilon_i $, 
the above problem can be formulated as  
\begin{subequations}
\label{e:OptProb-TM-baby}
\begin{align}
\max_{\tilde \varepsilon_i, \tilde \alpha, \tilde \beta, \theta} \quad &  2 \frac{\tilde \alpha - \tilde{\beta}}{\tilde \alpha + \tilde \beta}   \\
\textrm{subject to} \quad 
& \Pi^{i,\bkappa} \left( - \tilde \alpha  (\nabla + \imath \bkappa) \cdot (\nabla + \imath \bkappa)  - \tilde \varepsilon_i   \textrm{I} \right) \Pi^{i,\bkappa} \preceq  0, 
&& \bkappa \in \B_i, \  i \in [2] \\
\label{e:e:OptProb-TMc}
& (\textrm{I} - \Pi^{i,\bkappa}) \left( - \tilde \beta (\nabla + \imath \bkappa) \cdot (\nabla + \imath \bkappa)  - \tilde \varepsilon_i   \textrm{I} \right) ( \textrm{I} - \Pi^{i,\bkappa})  \succeq  0, 
&& \bkappa \in \B_i, \  i \in [2] \\
& \textrm{$\Pi^{i,\bkappa}$ is a rank $m$ orthogonal projection} 
&& \bkappa \in \B_i, \  i \in [2] \\
&  \theta \varepsilon_- \leq \tilde \varepsilon_i(\bx) \leq \theta \varepsilon_+, 
&& \bx \in \D, \ i \in [2] \\ 
& \theta \geq 0 \\ 
& (\varepsilon_1,\varepsilon_2) \in \mathcal{A}. 
%& \mbox{Idx}_m(\tilde \varepsilon_i / \theta ) = n_i, && i \in [2].
\end{align}
\end{subequations}
Finally, using homogeneity properties of the fractional objective function, we can equivalently rewrite \eqref{e:OptProb-TM-baby} with a linear objective function as follows (see, e.g., \cite[p. 151]{Boyd_2004}),   
\begin{subequations}
\label{e:OptProb-TM}
\begin{align}
\max_{\tilde \varepsilon_i, \tilde \alpha, \tilde \beta, \theta} \quad &  \tilde \alpha - \tilde \beta   \\
\textrm{subject to} \quad 
& \Pi^{i,\bkappa} \left( - \tilde \alpha (\nabla + \imath \bkappa) \cdot (\nabla + \imath \bkappa)  - \tilde \varepsilon_i   \textrm{I} \right) \Pi^{i,\bkappa} \preceq  0, 
&& \bkappa \in \B_i, \  i \in [2] \\
\label{e:e:OptProb-TMc2}
& (\textrm{I} - \Pi^{i,\bkappa}) \left( - \tilde \beta (\nabla + \imath \bkappa) \cdot (\nabla + \imath \bkappa)  - \tilde \varepsilon_i   \textrm{I} \right) ( \textrm{I} - \Pi^{i,\bkappa})  \succeq  0, 
&& \bkappa \in \B_i, \  i \in [2] \\
& \textrm{$\Pi^{i,\bkappa}$ is a rank $m$ orthogonal projection} 
&& \bkappa \in \B_i, \  i \in [2] \\
&  \theta \varepsilon_- \leq \tilde \varepsilon_i(\bx) \leq \theta \varepsilon_+, 
&& \bx \in \D, \ i \in [2] \\ 
& \tilde \alpha + \tilde \beta = 2, \ \ \theta \geq 0 \\ 
\label{e:OptProb-TM-Chern}
& (\varepsilon_1,\varepsilon_2) \in \mathcal{A}. 
%& \mbox{Idx}_m(\tilde \varepsilon_i / \theta ) = n_i, && i \in [2].
% & C_m^{(\bkappa_1, \bkappa_2)}(\tilde \varepsilon_i / \theta ) = n_i, && i \in [2]. 
\end{align}
\end{subequations}

We now approximate the solution to \eqref{e:OptProb-TM} with the following two-step method. Let $\Pi_m^{i,\bkappa}$ denote the projection
onto the eigenspace of the differential operator $A$ with $\bkappa \in \B$ and relative permittivity $\varepsilon_i$, $i \in [2]$, i.e., 
\begin{equation}
\label{e:ProjStep}    
\Pi_m^{i,\bkappa} \psi(\bx) = 
\sum_{j \in [m]} \left\langle  \phi_j^{(i)}(\bkappa; \cdot), \psi \right\rangle  \phi_j^{(i)}(\bkappa;\bx).
\end{equation}
With $\Pi = \Pi_m^{i,\bkappa} \psi(\bx)$ in \eqref{e:OptProb-TM} fixed, we solve 
\begin{subequations}
\label{e:OptProb-final}
\begin{align}
\max_{\tilde \varepsilon_i, \tilde \alpha, \tilde \beta, \theta} \quad &  \tilde \alpha - \tilde \beta   \\
\textrm{subject to} \quad 
& \Pi_m^{i,\bkappa} \left( - \tilde \alpha (\nabla + \imath \bkappa) \cdot (\nabla + \imath \bkappa)  - \tilde \varepsilon_i   \textrm{I} \right) \Pi_m^{i,\bkappa} \preceq  0, 
&& \bkappa \in \B_i, \  i \in [2] \\
\label{e:OptProb-final-C}
& (\textrm{I} - \Pi_m^{i,\bkappa}) \left( - \tilde \beta (\nabla + \imath \bkappa) \cdot (\nabla + \imath \bkappa)  - \tilde \varepsilon_i   \textrm{I} \right) ( \textrm{I} - \Pi_m^{i,\bkappa})  \succeq  0, 
&& \bkappa \in \B_i, \  i \in [2] \\
%& \textrm{$\Pi^{i,\bkappa}$ is a rank $m$ orthogonal projection} 
%&& \bkappa \in \B_i, \  i \in [2] \\
&  \theta \varepsilon_- \leq \tilde \varepsilon_i(\bx) \leq \theta \varepsilon_+, 
&& \bx \in \D, \ i \in [2] \\ 
& \tilde \alpha + \tilde \beta = 2, \ \ \theta \geq 0 \\ 
\label{e:OptProb-final-F}
& (\varepsilon_1,\varepsilon_2) \in \mathcal{A}. 
%& \mbox{Idx}_m(\tilde \varepsilon_i / \theta ) = n_i, && i \in [2].
% & C_m^{(\bkappa_1, \bkappa_2)}(\tilde \varepsilon_i / \theta ) = n_i, && i \in [2]. 
\end{align}
\end{subequations}
We then update $\Pi$ as in \eqref{e:ProjStep} and iterate. With the exception of the last constraint \eqref{e:OptProb-TM-Chern} on the topological invariant, this is a linear SDP in the variables $\tilde \alpha$, $\tilde \beta$, $\theta$, and $\tilde \varepsilon_i$.

%\subsection*{(TM) case} 
%Here we assume that $\rho_i \equiv 1$.    
% \begin{subequations}
% \label{e:OptProb-TM}
% \begin{align}
% \max_{\eta_i, \tilde{\alpha}, \tilde{\beta}} \quad &  \tilde{\alpha}^{-1} - \tilde{\beta}^{-1}   \\
% \textrm{subject to} \quad 
% & \Pi_m^{i,\bkappa} \left( - \tilde{\alpha}^{-1} \Delta  - \eta_i   \textrm{I} \right) \Pi_m^{i,\bkappa} \preceq  0, 
% && \bkappa \in \B_i, \  i \in [2] \\
% & (\textrm{I} - \Pi_m^{i,\bkappa}) \left( - \tilde{\beta}^{-1} \Delta  - \eta_i   \textrm{I} \right) ( \textrm{I} - \Pi_m^{i,\bkappa})  \succeq  0, 
% && \bkappa \in \B_i, \  i \in [2] \\
% & \textrm{$\Pi_m^{i,\bkappa}$ is a rank $m$ projection} 
% && \bkappa \in \B_i, \  i \in [2] \\
% & \tilde{\alpha}^{-1} + \tilde{\beta}^{-1} = 2  \\ 
% &  \eta_- \leq \eta_i \leq \eta_+, && i \in [2] \\
% & \textrm{Idx}(1,\eta_i) = n_i, && i \in [2]
% \end{align}
% \end{subequations}

To impose the constraint \eqref{e:OptProb-final-F}, we start with initial material coefficients that satisfy the constraint and use a proximal approach that penalizes the difference between the old and new iterations so that the topological invariant is preserved. We discuss the implementation of this strategy further for the topological indices considered in Section~\ref{s:TopologicalPhaseConstraint}.  Recall that the topological invariant is defined via the Berry curvature. Therefore, a linearization of the Berry curvature with respect to the electric permittivity $\varepsilon(\bx)$ is needed when the constraint is imposed in the new iteration from an old iteration. This is discussed in the following section.

\subsection{Sensitivity analysis of the Berry curvature}\label{section:SensitivityBerry}
To compute the Berry curvature $F_k(\bkappa)$ for the $k$-th band and its linearization, instead of using the definition \eqref{eq:Berry_curvature} directly, which involves the second derivative of the eigenmodes $\phi_j(\bkappa,\bx)$ with respect to the Bloch wave vector $\bkappa$ and is difficult to obtain computationally, we follow \cite{blanco2020tutorial} 
by combining the formulas \eqref{eq:Berry_phase1} and \eqref{eq:Berry_phase2} and considering the integral over a small region surrounding $\bkappa$.  More precisely, let $Q_{\bkappa}$ be a small square surround $\bkappa$ with the width $\delta \bkappa$. Then the integral \eqref{eq:Berry_phase1} and \eqref{eq:Berry_phase2} over $Q_{\bkappa}$ reads
\begin{equation*}
  \int_{Q_{\bkappa}}  F_k(\bkappa) \, d\bkappa = \oint_{\partial Q_{\bkappa}} \bzeta_k(\bkappa) \cdot d\bkappa,
\end{equation*}
which can be cast the following after discretization of the integral (\cite{blanco2020tutorial,vanderbilt2018berry}): 
\begin{equation}
\label{eq:F}
F_k(\bkappa) \cdot |Q_{\bkappa}|  = -\text{Im} \left(\ln \prod_{j=1}^4 \langle \phi(\bkappa_j;\cdot), \phi(\bkappa_{j+1};\cdot)\rangle\right)  + O(\delta\bkappa)^3.
\end{equation}
In the above, $\bkappa_j, \ j \in [4]$ represent four corners of the square $Q_{\bkappa}$ ordered counter-clockwisely and we set $\bkappa_5=\bkappa_1$. 
The above formula is gauge invariant in the sense that a phase change of $\phi(\bkappa,\bx)$ does not change the value of $F(\bkappa)$. Then the sensitivity analysis of the Berry curvature $F(\bkappa)$ is reduced to the sensitivity analysis of the eigenmode $\phi(\bkappa,\bx)$ over the corners of $Q_{\bkappa}$.

In the following, for clarity we neglect the high-order error $O(\delta\bkappa)^3$ arising from the discretization of the integral in \eqref{eq:F} in the computation of the Berry curvature and its sensitivity, and denote the approximate value by $\tF_m(\bkappa)$. That is, we set
\begin{equation}\label{eq:tF}
    \tF_k(\bkappa) = - \frac{1}{|Q_{\bkappa}|} \cdot \text{Im}\left(\ln \prod_{j=1}^4 \langle \phi(\bkappa_j;\cdot), \phi(\bkappa_{j+1};\cdot)\rangle\right).
\end{equation}

\begin{thm}
\label{s:SensitivityBerry}
For $\Lambda$-periodic permittivity $\varepsilon \in L^\infty(\mathbb{R}^2)$,
let  $\lambda(\bkappa) = \omega^2(\bkappa)$, $\bkappa \in \B$ 
be a dispersion surface with associated $\Lambda$-periodic eigenfunction $\phi(\bkappa; \cdot)$ satisfying \eqref{e:twistHelm}.  
Let  $\tF_k(\bkappa_0; \varepsilon)$ be the approximate Berry curvature of this surface at $\bkappa = \bkappa_0$ defined in \eqref{eq:tF}. 
Then $\tF_k(\bkappa_0; \varepsilon)$ is Fr\'echet differentiable with respect to $\varepsilon$ and there holds
\begin{equation}\label{eq:dtF}
\tF_k(\bkappa_0; \varepsilon + \delta \varepsilon) =
\tF_k(\bkappa_0; \varepsilon)  + \langle  \delta \varepsilon, g_k \rangle  
+ o\left( \| \delta \varepsilon \| \right), 
\end{equation}
where 
   \begin{equation}\label{eq:g}
       g_k(\bkappa_0; \bx) = \mbox{Im} \sum_{j=1}^{4}  \frac{1}{\langle \phi(\bkappa_j;\cdot),\phi(\bkappa_{j+1};\cdot)\rangle }
       \left(-\omega_{j}^{2}\phi(\bkappa_j;\bx)\overline{u_{j}(\bx)}-\omega_{j+1}^{2}v_{j}(\bx)\overline{\phi(\bkappa_{j+1};\bx)}+\phi(\bkappa_j;\bx)\overline{\phi(\bkappa_{j+1};\bx)}\right),
   \end{equation} 
$u_j$, $v_j \in H_{\varepsilon,p}^1(\D)$ solves the following equations in $\D$:
\begin{align*}
      & \frac{1}{\varepsilon} (\nabla+\imath\bkappa_{j})\cdot(\nabla+\imath\bkappa_{j}) u_j(\bx)+\omega_{j}^{2}u_j(\bx)+ 2\omega_j^2 \, \langle u_j, \phi(\bkappa_{j};\cdot) \rangle \, \phi(\bkappa_{j};\bx) = \phi(\bkappa_{j+1};\bx), \\
     & \frac{1}{\varepsilon} (\nabla+\imath\bkappa_{j+1})\cdot(\nabla+\imath\bkappa_{j+1}) v_j(\bx)+\omega_{j+1}^{2} v_j(\bx)+ 2\omega_{j+1}^2 \, \langle v_j, \phi(\bkappa_{j+1};\cdot) \rangle \, \phi(\bkappa_{j+1};\bx) = \phi(\bkappa_{j};\bx).
\end{align*}

\end{thm}

\begin{proof}
A direct perturbation of $\tF_k(\bkappa_0)$ in \eqref{eq:tF} yields
\begin{equation}\label{eq:pert_tF}
\delta \tF(\bkappa_0) \cdot |Q_{\bkappa}|  = -\text{Im}\sum_{j=1}^{4}
\frac{
\left(\langle \delta \phi(\bkappa_j;\cdot),\phi(\bkappa_{j+1};\cdot)\rangle +\langle \phi(\bkappa_j;\cdot),\delta \phi(\bkappa_{j+1};\cdot)\rangle +\int_{\D}\delta\varepsilon(\bx) \phi(\bkappa_j;\cdot)\overline{\phi(\bkappa_{j+1};\bx)} \, d\bx\right)}
{\langle \phi(\bkappa_j;\cdot),\phi(\bkappa_{j+1};\cdot)\rangle }.
\end{equation}
We compute $\langle \delta \phi(\bkappa_j;\cdot),\phi(\bkappa_{j+1};\cdot)\rangle $ and $\langle \phi(\bkappa_j;\cdot),\delta \phi(\bkappa_{j+1};\cdot)\rangle $ by using the adjoint method as follows.

First, a perturbation of the equation \eqref{e:twistHelm} at $\bkappa=\bkappa_j$ with respect to $\varepsilon$ leads to 
\begin{equation}
(\nabla+i\bkappa_{j})\cdot(\nabla+i\bkappa_{j})\delta \phi(\bkappa_j;\bx)+\omega_{j}^{2}\varepsilon(\bx)\delta \phi(\bkappa_j;\bx)=-2\omega_{j}\delta\omega_{j}\varepsilon(\bx)\phi(\bkappa_j;\bx)-\omega_{j}^{2}\delta\varepsilon(\bx)\phi(\bkappa_j;\bx).\label{eq:perturb equation}
\end{equation}
By choosing $\norm{\phi(\bkappa_j;\cdot)}_{L^2_{\varepsilon,p}(\D)}=1$, it follows from \eqref{e:twistHelm}, \eqref{eq:perturb equation} and the Green's formula that the perturbation of the frequency is given by
\begin{equation}\label{eq:delta_omega}
\delta\omega_{j}=-\frac{\omega_{j}}{2}\int_{\D}\delta\varepsilon(\bx)|\phi(\bkappa_j;\bx)|^{2} \, d\bx. 
\end{equation}

%\[
% \int\delta\varepsilon|\phi(\bkappa_j;\cdot)|^{2}dx+\langle \delta \phi(\bkappa_j;\cdot),\phi(\bkappa_j;\cdot)\rangle +\langle \phi(\bkappa_j;\cdot),\delta \phi(\bkappa_j;\cdot)\rangle =0.
% \]
% \[
% \int\delta\varepsilon|\phi(\bkappa_j;\cdot)|^{2}dx+\langle \delta \phi(\bkappa_j;\cdot),\phi(\bkappa_j;\cdot)\rangle +\overline{\langle \delta \phi(\bkappa_j;\cdot),\phi(\bkappa_j;\cdot)\rangle }=0.
% \]

For clarity let us denote the self-adjoint differential operator $\frac{1}{\varepsilon} (\nabla+\imath\bkappa_{j})\cdot(\nabla+\imath\bkappa_{j})+\omega_{j}^{2}$ by $\Lo_j$  and rewrite \eqref{eq:perturb equation} as
\begin{equation}\label{eq:pert_phi}
    \Lo_j \delta \phi(\bkappa_j;\bx) = -2\omega_{j}\delta\omega_{j}\phi(\bkappa_j;\bx)-\omega_{j}^{2} \frac{\delta\varepsilon(\bx)}{\varepsilon(\bx)}\phi(\bkappa_j;\bx).
\end{equation}
We introduce the adjoint equation
\begin{equation}\label{eq:adjoint_phi1}
    \Lo_j u_j(\bx) = \phi(\bkappa_{j+1};\bx) - \langle \phi(\bkappa_{j+1};\cdot), \phi(\bkappa_{j};\cdot) \rangle \, \phi(\bkappa_{j};\bx),
\end{equation}
in which the right side is orthogonal to $\mbox{Ker}(\Lo_j)$. By the Fredholm alternative, \eqref{eq:adjoint_phi1} 
attains a solution in $H_{\varepsilon, p}^{1}(\D)$. Alternatively, we choose a solution $u_j\in H_{\varepsilon, p}^{1}(\D)$ that satisfies
\begin{equation}\label{eq:adjoint_phi2}
  \tilde \Lo_j u_j(\bx) :=  \Lo_j v(\bx) + t \, \langle u_j, \phi(\bkappa_{j};\cdot) \rangle \, \phi(\bkappa_{j};\bx) = \phi(\bkappa_{j+1};\bx),
\end{equation}
wherein $t\neq0$ is any fixed real number. It can be shown that \eqref{eq:adjoint_phi2} attains a unique solution in $H_{\varepsilon,p}^{1}(\D)$. In addition, in view of \eqref{eq:adjoint_phi1} and \eqref{eq:adjoint_phi2}, it is clear that
\begin{equation}\label{eq:u_proj}
    t \, \langle u_j, \phi(\bkappa_{j};\cdot) \rangle = \langle \phi(\bkappa_{j+1};\cdot), \phi(\bkappa_{j};\cdot) \rangle.
\end{equation}
Although \eqref{eq:adjoint_phi1} and \eqref{eq:adjoint_phi2} are equivalent when \eqref{eq:u_proj} holds, \eqref{eq:adjoint_phi2} is computationally more convenient as the kernel of the operator $\tilde \Lo_j: H_{\varepsilon,p}^{1}(\D) \to H^{-1}(\D)$ is $span\{0\}$, while
$\Lo_j: H_{\varepsilon,p}^{1}(\D) \to H^{-1}(\D)$ attains a nontrivial kernel space.

Using \eqref{eq:pert_phi} and \eqref{eq:adjoint_phi2}, we have 
\begin{align}\label{eq:green1}
     & \langle \delta\phi(\bkappa_j;\cdot), \tilde \Lo_j u_j \rangle - \langle \Lo_j\delta\phi(\bkappa_j;\cdot), u_j \rangle \\
     = & \, \langle \delta \phi(\bkappa_j;\cdot),\phi(\bkappa_{j+1};\cdot)\rangle 
     + 2 \delta\omega_j \omega_j \langle \phi(\bkappa_j;\cdot), u_j \rangle + \omega_j^2 \int_{\D} \delta\varepsilon(\bx) \phi(\bkappa_j;\cdot) \overline{u_j (\bx)}  \, d\bx. \nonumber
\end{align}
On the other hand, it follows from the Green's formula that
\begin{equation}\label{eq:green2}
 \langle \delta\phi(\bkappa_j;\cdot), \tilde \Lo_j u_j \rangle - \langle \Lo_j\delta\phi(\bkappa_j;\cdot), u_j \rangle = t \, \langle \phi(\bkappa_{j};\cdot), u_j \rangle \langle \delta \phi(\bkappa_j;\cdot), \phi(\bkappa_{j};\cdot)\rangle.
\end{equation}
Therefore, with a combination of \eqref{eq:green1}, \eqref{eq:green2} and \eqref{eq:delta_omega}, we obtain
\begin{align}\label{eq:dphi_phi1}
    & \langle \delta \phi(\bkappa_j;\cdot),\phi(\bkappa_{j+1};\cdot)\rangle \\
    = & \, \omega_j^2 \langle \phi(\bkappa_j;\cdot), u_j \rangle \int_{\D}\delta\varepsilon(\bx)|\phi(\bkappa_j;\bx)|^{2} \, d\bx
    + t \, \langle \phi(\bkappa_{j};\cdot), u_j \rangle \langle \delta \phi(\bkappa_j;\cdot), \phi(\bkappa_{j};\cdot)\rangle - \omega_j^2 \int_{\D} \delta\varepsilon(\bx) \phi(\bkappa_j;\cdot) \overline{u_j (\bx)}  \, d\bx \nonumber. 
\end{align}
Note that $\norm{\phi(\bkappa_j;\cdot)}_{L_{\varepsilon,p}^2(\D)}=1$ and $\norm{\phi(\bkappa_j;\cdot)+\delta\phi(\bkappa_j;\cdot)}_{L_{\varepsilon,p}^2(\D)}=1$, it can be deduced that
\[
\langle \delta \phi(\bkappa_j;\cdot),\phi(\bkappa_j;\cdot)\rangle =-\frac{1}{2}\int_{\D}\delta\varepsilon|\phi(\bkappa_j;\cdot)|^{2}dx+\imath\alpha_j
\]
for some real number $\alpha_j$. By substituting the above into \eqref{eq:dphi_phi1} and choosing $t=2\omega_j^2$, we arrive at
\begin{align}\label{eq:dphi_phi2}
    \langle \delta \phi(\bkappa_j;\cdot),\phi(\bkappa_{j+1};\cdot)\rangle 
    & = - \omega_j^2 \int_{\D} \delta\varepsilon(\bx) \phi(\bkappa_j;\cdot) \overline{u_j (\bx) } \, d\bx +  \imath \alpha_j t \langle \phi(\bkappa_j;\cdot), u_j \rangle \nonumber \\
    & =  - \omega_j^2 \int_{\D} \delta\varepsilon(\bx) \phi(\bkappa_j;\cdot) \overline{u_j (\bx) } \, d\bx +  \imath \alpha_j \langle \phi(\bkappa_{j};\cdot), \phi(\bkappa_{j+1};\cdot) \rangle,
\end{align}
where we have used the relation \eqref{eq:u_proj}.

Following parallel lines, it can be shown that
\begin{equation}\label{eq:phi_dphi}
    \langle \phi(\bkappa_j;\cdot),\delta \phi(\bkappa_{j+1};\cdot)\rangle 
    = - \omega_{j+1}^2 \int_{\D} \delta\varepsilon(\bx)  v_j(\bx) \overline{\phi(\bkappa_{j+1};\cdot)}  \, d\bx  -  \imath\alpha_{j+1} \langle \phi(\bkappa_{j};\cdot), \phi(\bkappa_{j+1};\cdot) \rangle,
\end{equation}
where $\alpha_{j+1}$ is some real number and $v_j\in H_{\varepsilon,p}^1(\D)$ solves
\begin{equation*}
  \tilde \Lo_{j+1} v_j(\bx) :=  \Lo_{j+1} v_j(\bx) +  2\, \omega_{j+1}^2 \langle v_j, \phi(\bkappa_{j+1};\cdot) \rangle \, \phi(\bkappa_{j+1};\bx) = \phi(\bkappa_{j};\bx).
\end{equation*}
Substituting \eqref{eq:dphi_phi2} and \eqref{eq:phi_dphi} into \eqref{eq:pert_tF} completes the proof.
\end{proof}

\subsection{Numerical Implementation} 
\label{s:NumImp}
\subsubsection{Forward PDE solver} \label{ss: PDE solver} For a given permittivity function $\varepsilon(\bx)$, the eigenvalue problem 
 (\ref{e:twistHelm}a-\ref{e:twistHelm}b) can be efficiently solved by finite difference methods \cite{kao2005maximizing,kao2019extremal} or finite element methods \cite{dobson1999efficient,dobson2000efficient}. Here we adopt the approach in \cite{kao2019extremal} which first transforms the fundamental periodic cell $\D$ to a unit square by a linear mapping and then apply the standard nine-point finite difference discretization of the transformed equation on the unit square to achieve second-order accurate solutions.

\subsubsection{Topological phase constraint} \label{s:TopologicalPhaseConstraint}
First, the computation of Berry curvature $F_k(\bkappa)$ is performed by taking the leading-order term in \eqref{eq:F} which requires the eigenmodes for different wave vectors. To enforce the constraint \eqref{e:OptProb-TM-Chern}, we focus the discussion on the case
when the topological index of the photonic crystal is given by the gap valley Chern number with $\mbox{Idx}_m= C_{g,m}^{(\bkappa_1, \bkappa_2)}$, defined in \eqref{eq:gap_valley_Chern_number}. Recall that the gap Berry curvature $F_{g,m}(\bkappa)=\sum_{\ell\in[m]} F_\ell(\bkappa)$ and $F_\ell(\bkappa)$ is the Berry curvature for the $\ell$-th band,
we may replace the constraint \eqref{e:OptProb-TM-Chern}  by 
\begin{equation}\label{eq:constraint_F}
  \sgn\big(F_{g,m}(\bkappa_j;\tilde\varepsilon_i)\big) \cdot F_{g,m}(\bkappa_j;\tilde\varepsilon_i) \ge \tau_0 
  \sgn\big(F_{g,m}(\bkappa_j;\tilde\varepsilon_i^{[0]}) \big) \cdot F_{g,m}(\bkappa_j;\tilde\varepsilon_i^{[0]}),
  \quad i,j\in [2];
\end{equation}
for some chosen real number $0<\tau_0\le1$, wherein $\tilde\varepsilon_i^{[0]}$ represents the initial guess for $\tilde\varepsilon_i$.  The constraint above will keep $|F_{g,m}(\bkappa_j;\tilde\varepsilon_i)|$ strictly positive during the iteration and will reserve valley Chern number for the optimized periodic media. At each iteration $p$, the Berry curvature $F_{g,m}(\bkappa)$ is linearized and \eqref{eq:constraint_F} is replaced by 
\begin{equation}\label{eq:constraint_F_linear}
  \sgn\big(F_{g,m}(\bkappa_j;\tilde\varepsilon_i^{[p]})\big) \cdot \delta F_{g,m}(\bkappa_j;\tilde\varepsilon_i^{[p]}) + \tau_p
  \sgn\big(F_{g,m}(\bkappa_j;\tilde\varepsilon_i^{[0]}) \big) \cdot F_{g,m}(\bkappa_j;\tilde\varepsilon_i^{[0]}) \ge 0,
  \quad i,j\in [2];
\end{equation}
for $\tau_p$ satisfying $\sum_{p=0}^\infty \tau_p=1-\tau_0$. In the practical implementation, $\delta F_{g,m}=\sum_{\ell\in[m]} \delta F_\ell(\bkappa)$, wherein the perturbation $\delta F_\ell(\bkappa)$ is replaced by its leading-order given in \eqref{eq:dtF}, and we choose a geometric series $\{\tau_p\}_{p=0}^\infty$.

\medskip
\subsubsection{Optimization algorithm}
In summary, the semi-definite optimization approach for solving the problem \eqref{e:OptProb} is summarized in Algorithm~\ref{alg:sem_define_opt}. We use the \texttt{CVX} optimization package \cite{grant2014cvx} which is a MATLAB-based modeling system to solve \eqref{e:OptProb-final}. In our implementation of the constraint \eqref{e:OptProb-final-C}, the spectral projection $\textrm{I} - \Pi_m^{i,\bkappa}$ onto the space associated to the dispersion surfaces above $\lambda_{m+1}^{(i)}(\bkappa)$ is truncated using a fixed number of surfaces. 
Generally, we find that the number of surfaces that should be used matches the largest multiplicity eigenvalue in the dispersion surface $\lambda_{m+1}^{(i)}(\bkappa)$. For our implementation, we use $3$ surfaces.  
To respect the topological constraints \eqref{e:OptProb-TM-Chern} numerically, the linear inequality constraint (\ref{eq:constraint_F_linear}) is implemented. The convergence (stopping) criterion is based on that the norm of two consecutive differences for $\varepsilon_1$ and $\varepsilon_2$ are both less than given tolerance.

\begin{algorithm}[t!]
\label{alg:sem_define_opt}
	\caption{Maximizing the shared gap-to-midgap ratio $G_m(\varepsilon_1,\varepsilon_2) $ }
	\begin{algorithmic}
        \State Initial guess: choose an admissible pair of permittivities $(\varepsilon_1,\varepsilon_2) \in \mathcal{A}$ (defined in Section~\ref{ss:forOptim}).
        \While {not converged}
        \begin{enumerate}
        \item  For the current $(\varepsilon_1,\varepsilon_2) \in \mathcal{A}$, compute their corresponding eigenvalues and eigenfunctions, satisfying
 \eqref{e:twistHelmA} -- \eqref{e:twistHelmB} on the fundamental periodic cell $\D$ for specified $\bkappa \in \B$  via a finite difference method discussed in Section~\ref{ss: PDE solver}.
          \item Given $m$ and $\bkappa \in \B$, compute the projection $\Pi_m^{i,\bkappa}$  for $i \in [2]$ by \eqref{e:ProjStep}.  
          \item Compute the gap Berry curvature $F_{g,m}(\bkappa;\varepsilon_i)$ and its linearization  for $i \in [2]$ using \eqref{eq:F} and \eqref{eq:dtF}.
          \item
          Solve \eqref{e:OptProb-final} to find a new pair of permittivities $(\varepsilon_1,\varepsilon_2) \in \mathcal{A}$
          while respecting the linearized topological constraint \eqref{eq:constraint_F_linear}.
          \item Update the permittivities $(\varepsilon_1,\varepsilon_2) \in \mathcal{A}$. 
        \end{enumerate}
		\EndWhile
	\end{algorithmic}
\end{algorithm}

\subsubsection{Topological edge mode}
The edge mode for a joint structure formed by two photonic crystals satisfies the spectral problem \eqref{eq:edge_mode_PDE} -- \eqref{eq:edge_mode_decay}. In light of the quasi-periodic condition \eqref{eq:edge_mode_qp_bnd},  the problem is reduced to an infinite strip  $D_a:=\{\bx; \; 0<\bx\cdot\ba<a \}$, in which $a$ is the lattice constant along the interface direction $\ba$. We perform the numerical computation of edge modes via the super cell method by truncating the infinite strip into a finite large one given by $D_{a,L}:=\{\bx; \; 0<\bx\cdot\ba<a  , |\bx \cdot \ba^\perp| < L \}$ and imposing the zero Dirichlet boundary condition on the side boundary of the strip. Due to the exponential decay of the edge modes, it can be shown that the solution of the spectral problem in $D_{a,L}$ converges to that of the spectral problem \eqref{eq:edge_mode_PDE} -- \eqref{eq:edge_mode_decay} with the same rate as $L\to\infty$ \cite{soussi2005convergence}.

\section{Computational results} 
\label{s:CompResults}
\subsection{Square lattice} \label{s:ex1}
\begin{figure}[!p]
\begin{centering}
\includegraphics[scale=1,trim = 0mm 5mm 0mm 5mm, clip]{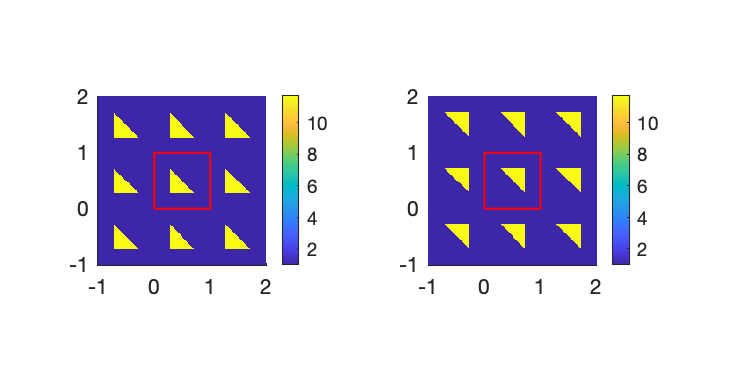}
\includegraphics[scale=0.65,trim = 0mm 5mm 0mm 5mm, clip]{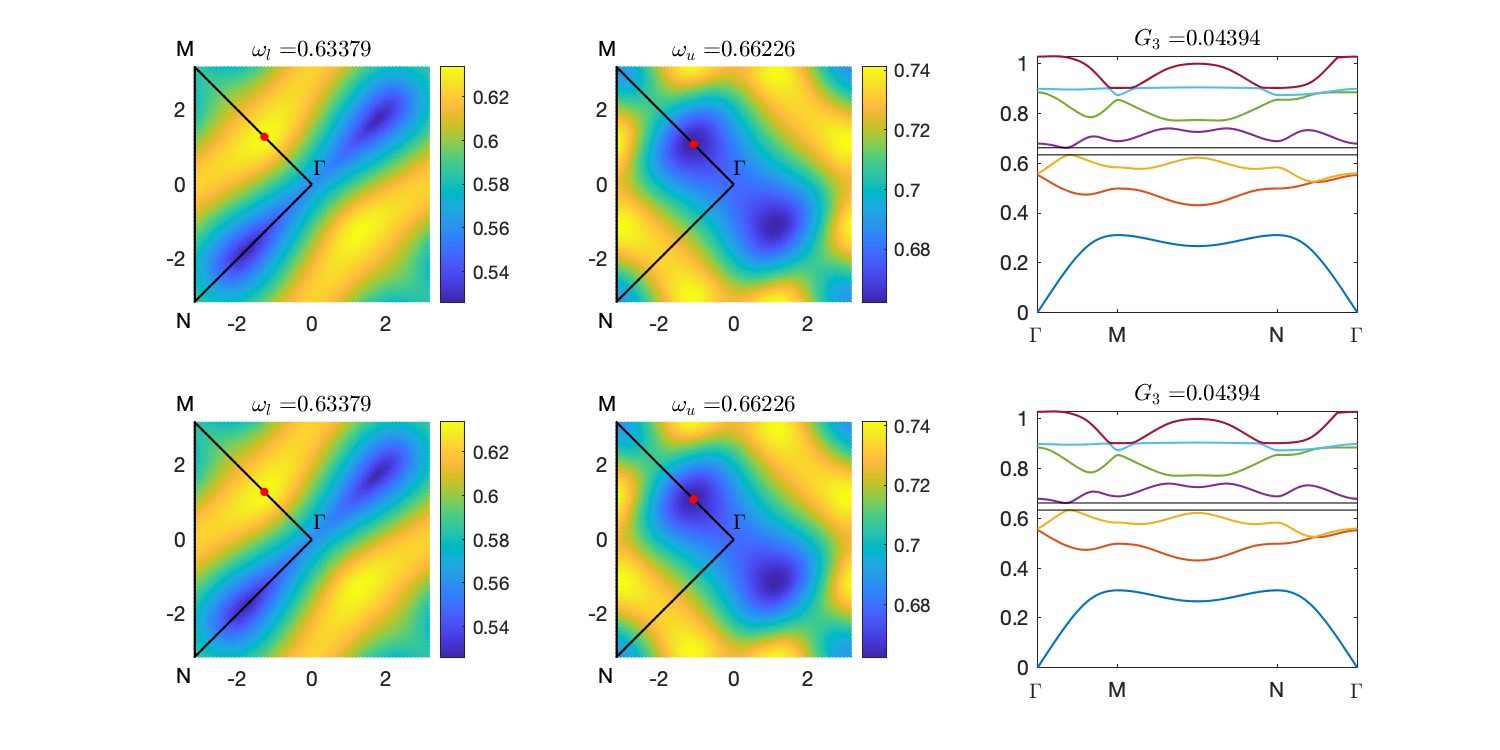}
\includegraphics[width=5.5cm,,trim = 0mm 5mm 0mm 5mm, clip]{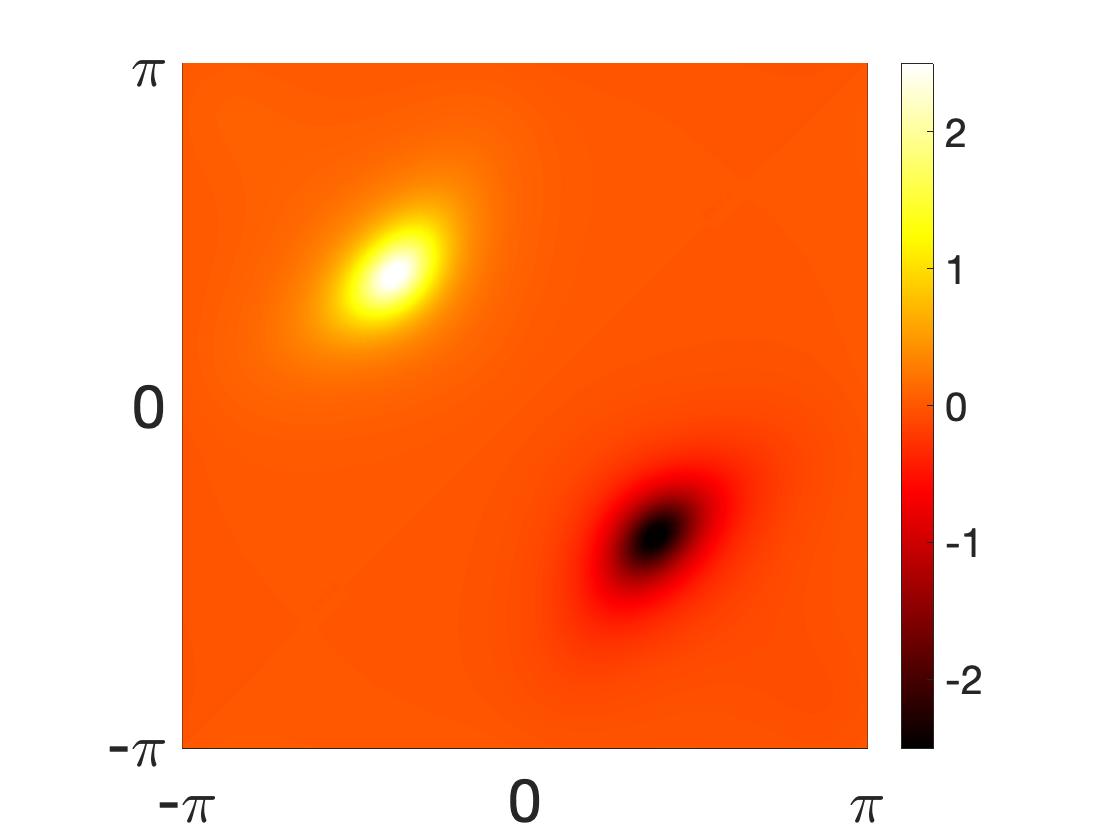}
\includegraphics[width=5.5cm,trim = 0mm 5mm 0mm 5mm, clip]{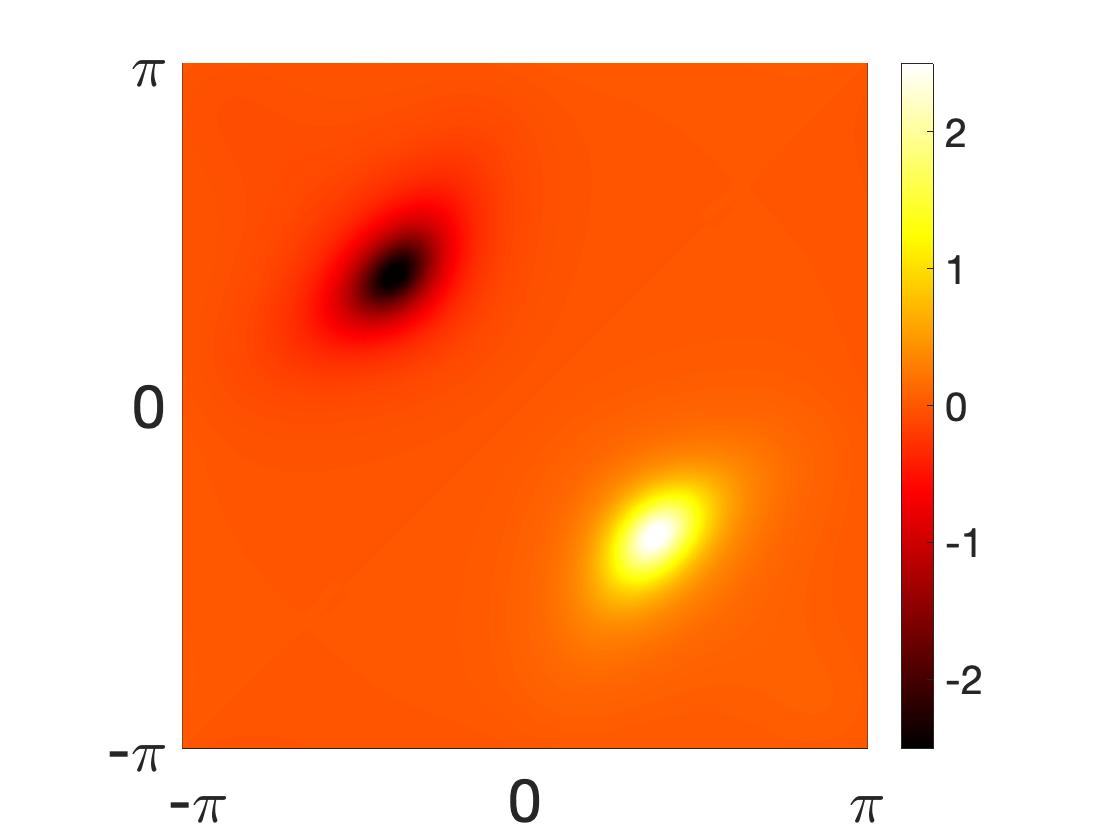}
\par\end{centering}
\caption{Row 1: Distribution of the electric permittivity for the initial configuration of photonic crystals.
Row 2 and 3: The third (left column) and fourth (middle column) band of the two photonic crystals. The eigenvalues along the boundary of the reduced Brillouin zone $\Gamma$MN for two photonic crystals is plotted on the right column.
Row 4: The gap Berry curvature $F_{g,3}(\bkappa)$ for two initial photonic crystals. The peak values of $F_{g,3}(\bkappa)$ appear at $\bkappa_1\approx(-1.2, 1.2)$ and $\bkappa_2\approx(1.2, -1.2)$. See Section~\ref{s:ex1}.
}
\label{Fig:square_lattice_init_PCs}
\end{figure}

In this example, we consider photonic crystals with the periodicity of a square lattice, for which the lattice vectors are given by $\be_1=(1,0)$ and $\be_2=(0,1)$. 
Figure \ref{Fig:square_lattice_init_PCs} shows two initial PCs. 
Each PC consists of an array of inclusions with geometry given by an isosceles right triangle whose short edge length is 0.45  \cite{kim2021multiband}. 
The orientation of the inclusions is different for the two PCs: in the first PC, the right angle of the triangle is at the lower left while in the second PC it is at the upper right.  
Here and hereafter, the permittivity value of the inclusions is assumed to be $\varepsilon_+ = 11.7$ while the background permittivity is $\varepsilon_- = 1$. 

In this example, we will maximize the gap between the third and fourth band of the spectrum, which is opened through the perturbation of a Dirac point by breaking the symmetry of inclusions \cite{kim2021multiband}.
The  third and fourth dispersion surfaces corresponding to these initial PC are plotted in 
Figure \ref{Fig:square_lattice_init_PCs} (Row 2-3).
It can be calculated that the initial shared gap-to-midgap is $G_{3}\approx 0.04394$. 
The corresponding gap Berry curvature $F_{g,3}(\bkappa)$ for the two PCs is plotted in Figure \ref{Fig:square_lattice_init_PCs} (Row 4), where the peak values appear at $\bkappa_1\approx(-1.2, 1.2)$ and $\bkappa_2\approx(1.2, -1.2)$.  The valley Chern number, which is the topological invariant, takes the value $1$ and $-1$ for the two PCs respectively.

\begin{figure}[!p]
\begin{centering}
\includegraphics[scale=1,trim = 0mm 5mm 0mm 5mm, clip]{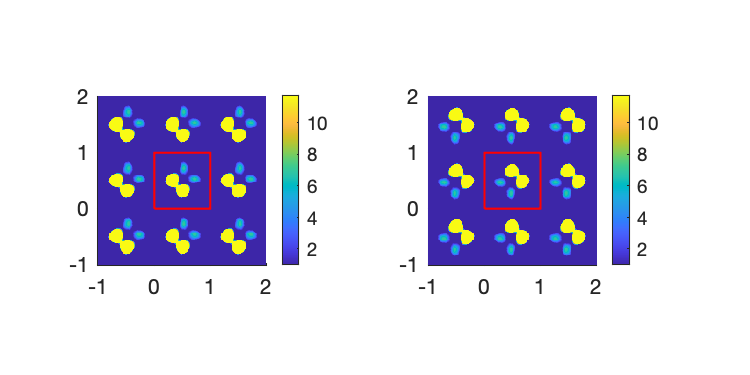} 
\includegraphics[scale=0.65,trim = 0mm 5mm 0mm 5mm, clip]{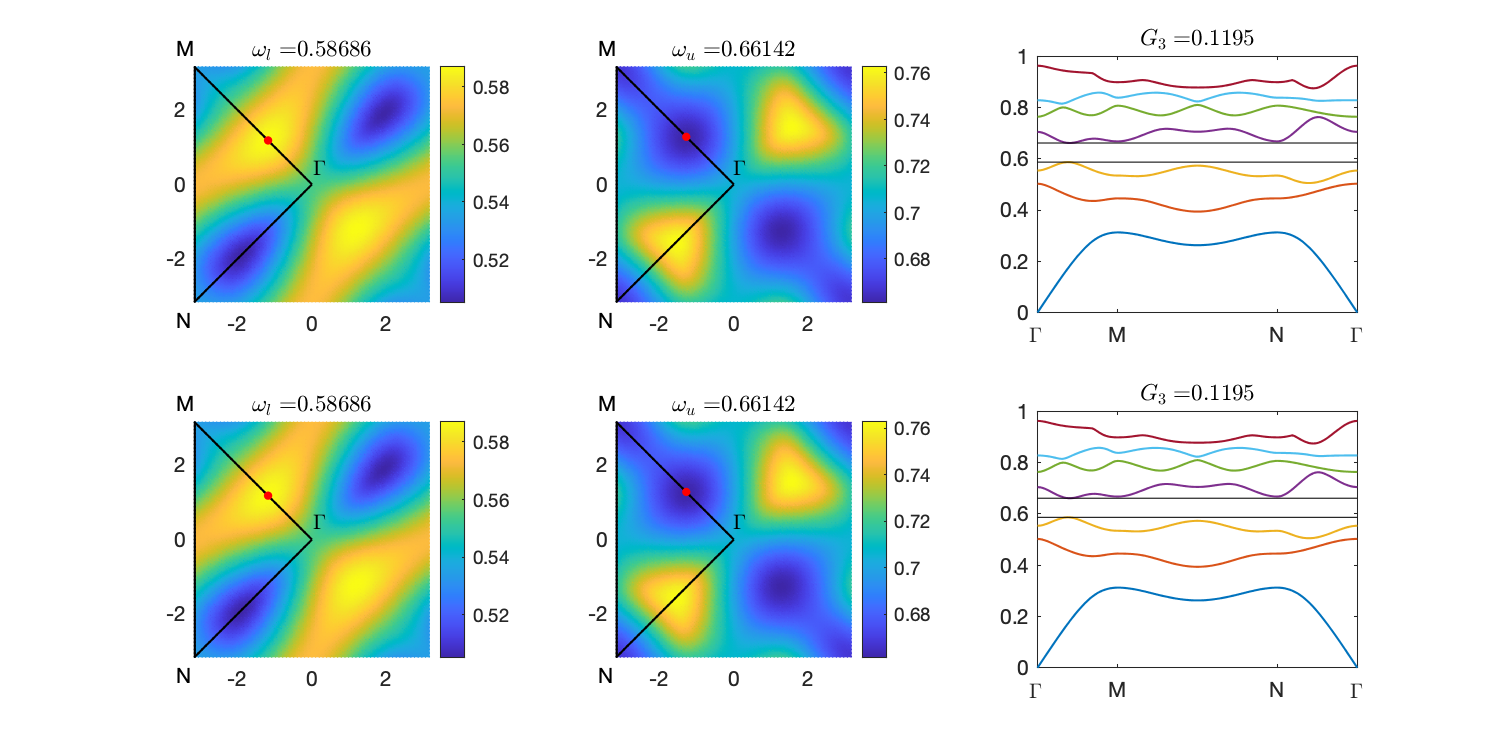}
\includegraphics[width=5.5cm,trim = 0mm 5mm 0mm 5mm, clip]{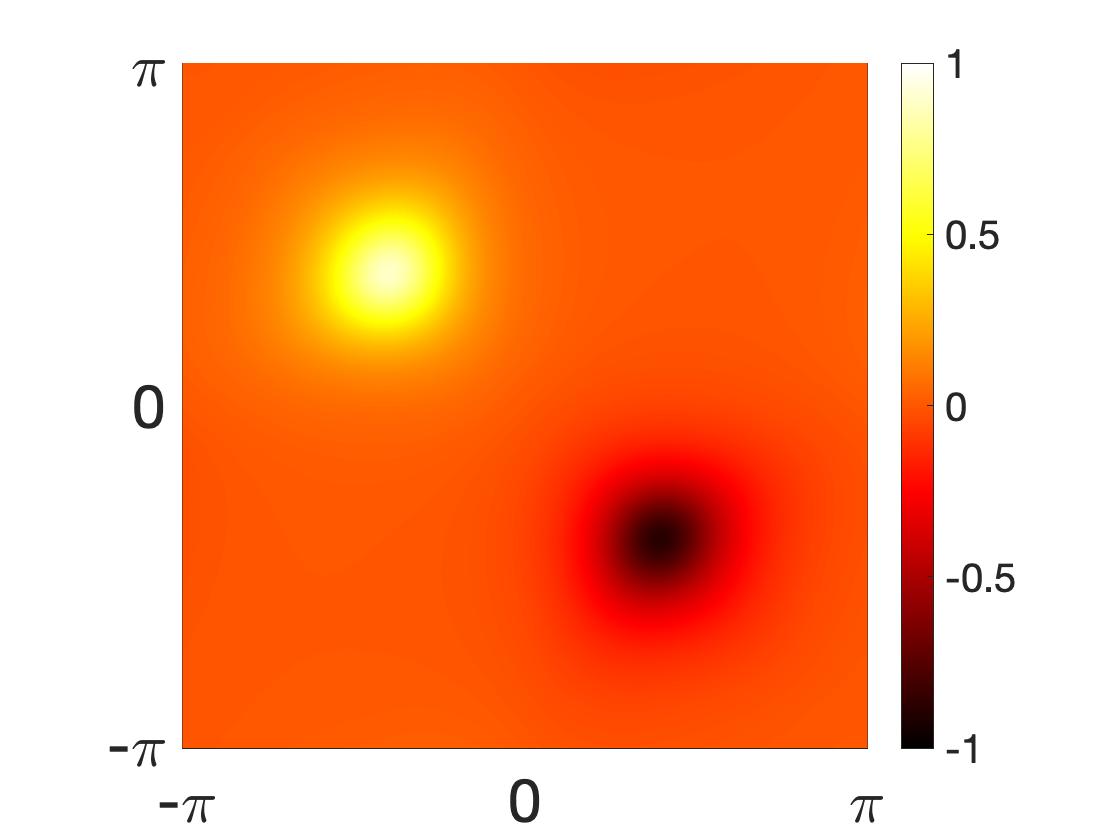}
\hspace*{-20pt}
\includegraphics[width=5.5cm,trim = 0mm 5mm 0mm 5mm, clip]{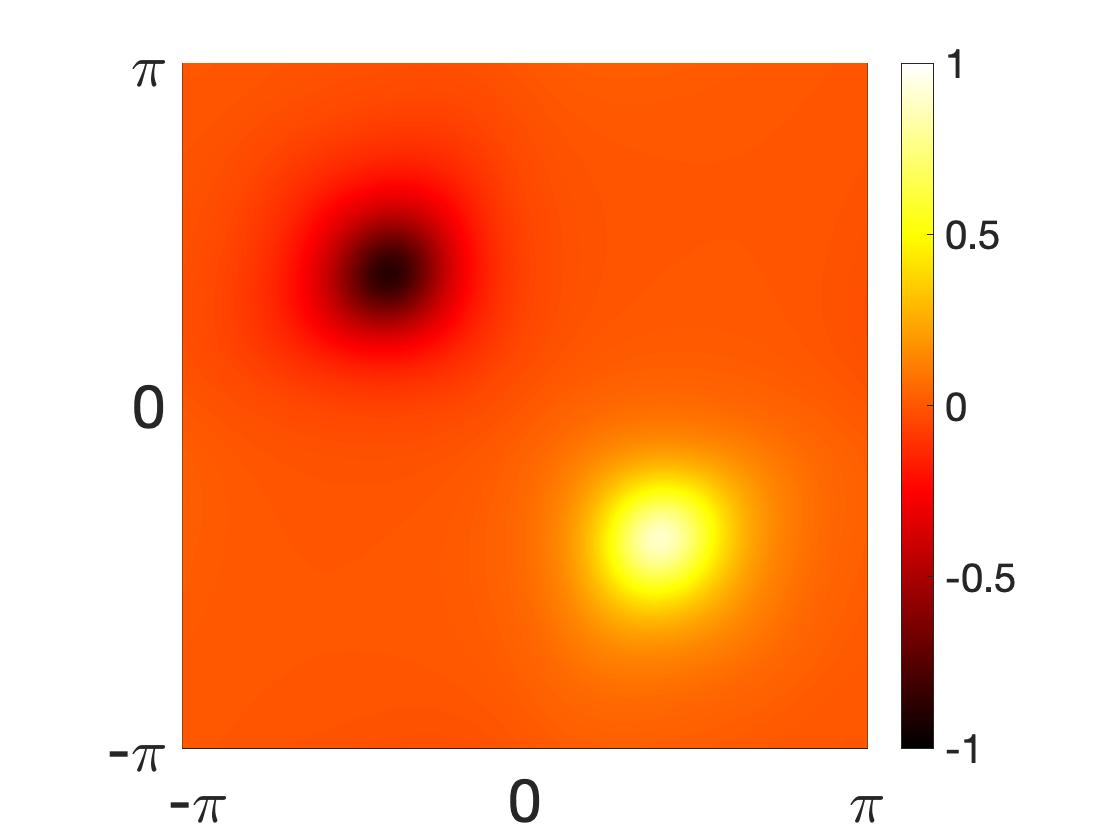}
\par\end{centering}
\caption{Row 1: Distribution of the electric permittivity for the optimized photonic crystals and the corresponding band structures.
Row 2 and 3: The third (left column) and fourth (middle column) band of two optimized photonic crystals. The eigenvalues along the boundary of the reduced Brillouin zone $\Gamma$MN for two photonic crystals is plotted on the right column.
Row 4: The gap Berry curvature $F_{g,3}(\bkappa)$ for two optimized photonic crystals.
The peak values of Berry curvature appear at $\bkappa_1\approx(-1.2, 1.2)$ and $\bkappa_2\approx(1.2, -1.2)$; see Section~\ref{s:ex1}.
}
\label{Fig:square_lattice_opt_PCs}
\end{figure}

Our implementation of Algorithm 1 terminates in 10 iterations leading to the optimized PCs shown in Figure \ref{Fig:square_lattice_opt_PCs}. 
The inclusions for the optimized PCs consist of four ovals in each periodic cell. Note also that intermediate values between $\varepsilon_-$ and $\varepsilon_+$ are attained. 
In the implementation, we set $\tau_0=\frac{1}{2}$ in \eqref{eq:constraint_F} and found that the band gap does not  significantly improve when $\tau_0$ is reduced unless it is close to 0, for which the gap Berry curvature essentially  becomes trivial. 
The band structure of the optimized PCs is shown in Figure \ref{Fig:square_lattice_opt_PCs} (Row 2-3), and the shared gap-to-midgap ratio is $G_{3}\approx 0.1195$. The gap Berry curvature $F_{g,3}(\bkappa)$ for two PCs is plotted in Figure \ref{Fig:square_lattice_opt_PCs} (Row 4). It is clear that the valley Chern number is preserved during the iterations. 

\begin{figure}[!p]
\begin{centering}
\includegraphics[width=16cm]{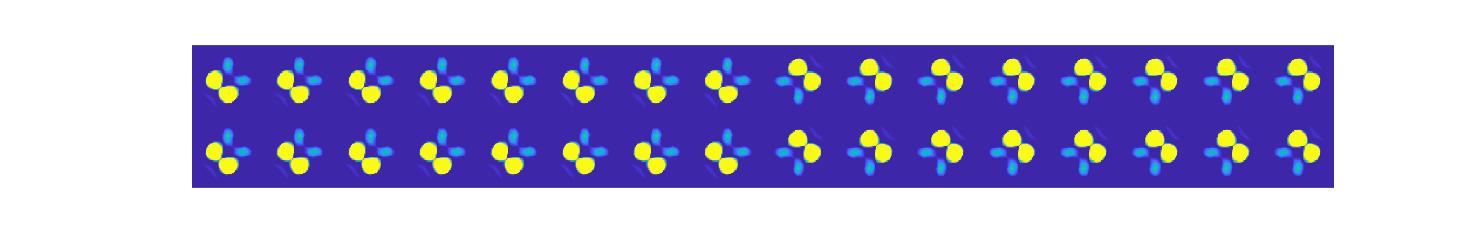}
\includegraphics[width=8cm]{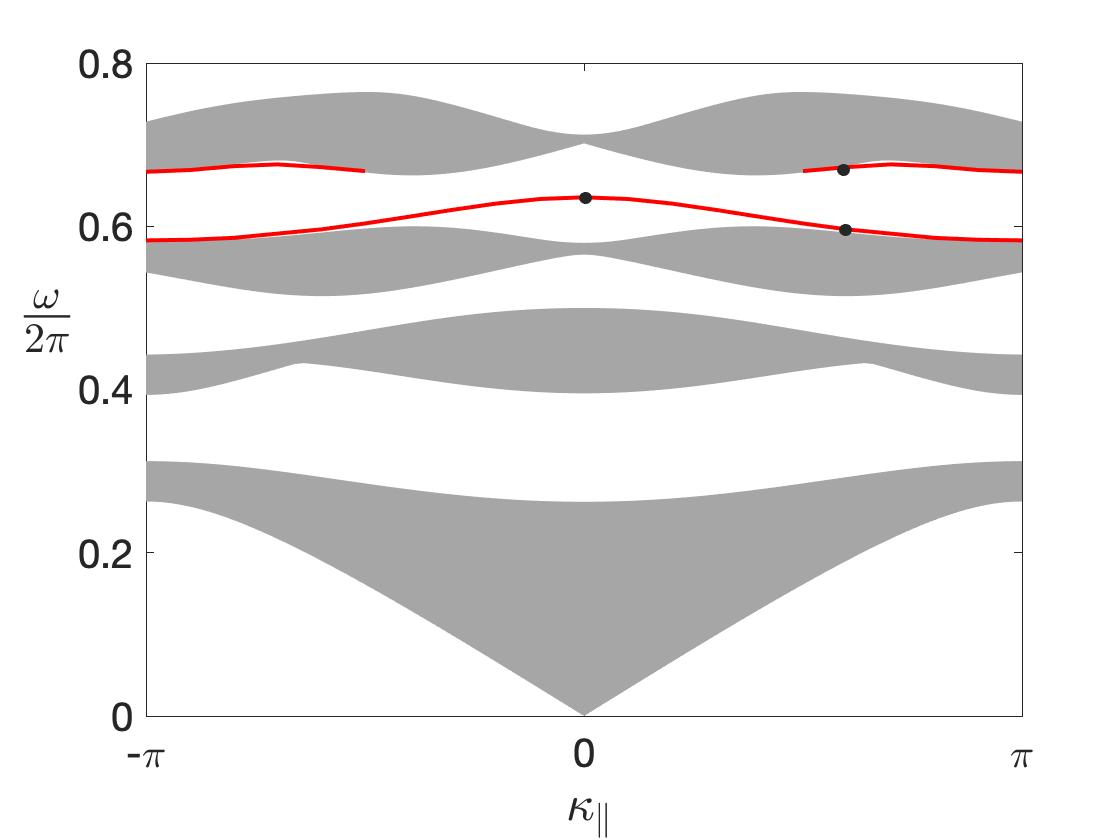}
\includegraphics[width=16cm]{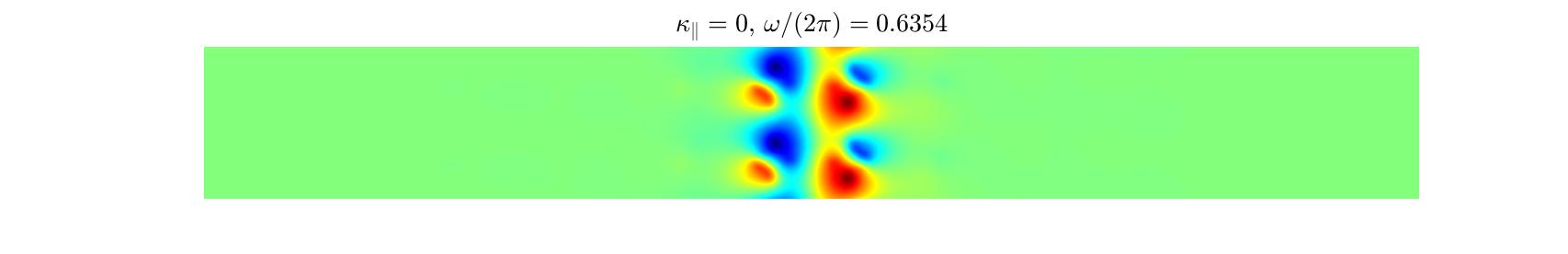}
\includegraphics[width=16cm]{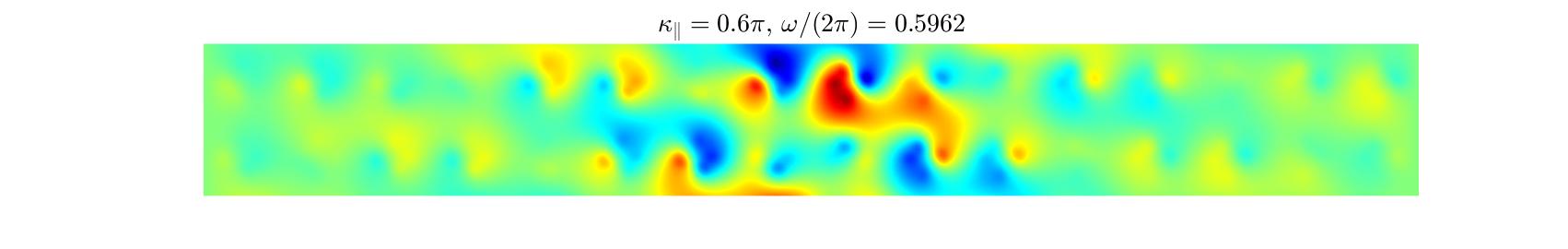}
\includegraphics[width=16cm]{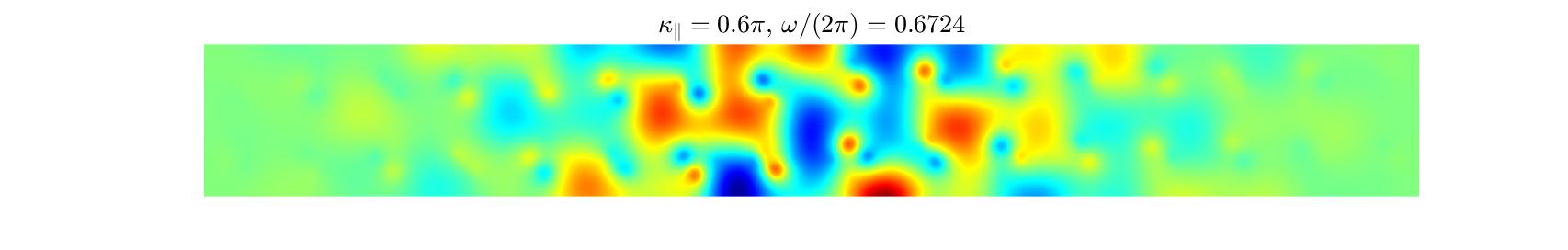}
\par\end{centering}
\caption{Row 1: The joint structure formed by gluing two optimized periodic media together along the $x_2$-axis. Row 2: The dispersion relation of the edge modes along the interface for the wavenumber $\kappa_\parallel\in[-\pi,\pi]$. The gray area is the projection of the continuous spectrum of the periodic media along the interface direction.
Row 3-5: Real part of the edge modes when $\kappa_\parallel=0$ and $0.6\pi$ respectively. The corresponding eigenfrequencies are marked as black dots in the dispersion relation. 
See Section~\ref{s:ex1}.} 
\label{Fig:square_lattice_joint_PC}
\end{figure}

Finally,  we consider the joint structure by gluing the two optimized periodic media along the $x_2$-axis as shown in Figure \ref{Fig:square_lattice_joint_PC} (Row 1). Due to the topological phase transition across the interface, edge modes can be supported along the interface of the joint structure. 
Figure \ref{Fig:square_lattice_joint_PC} (Row 2) plots the dispersion relations for the edge modes, where $\kappa_\parallel\in[-\pi, \pi]$ is the wavenumber along the interface direction. It is observed the spectral bandwidth of edge modes spans almost the whole spectral gap, thus it is significantly enlarged after the optimization. We also plot
the edge modes when $\kappa_\parallel=0$ and $0.6\pi$ in Figure \ref{Fig:square_lattice_joint_PC} (Row 3-5). These eigenmodes are localized at the interface of the joint structure and decay exponentially along the positive and negative $x_1$-axis.

% \begin{figure}[!p]
% \begin{centering}
% \includegraphics[scale=0.8,trim = 0mm 5mm 0mm 5mm, clip]{figs/square_ini_V.png}
% \par\end{centering}
% \vspace*{-20pt}
% \caption{Distribution of the electric permittivity for the initial configuration of photonic crystals; see Section~\ref{s:ex1}.}
% \label{Fig:square_lattice_init_PCs}
% \end{figure}

% \begin{figure}[!htbp]
% \begin{centering}
% \includegraphics[scale=0.65,trim = 0mm 5mm 0mm 5mm, clip]{figs/square_ini_bandgap.png}
% \par\end{centering}
% \caption{The third (left column) and fourth (middle column) band of the two photonic crystals in Figure \ref{Fig:square_lattice_init_PCs}. The eigenvalues along the boundary of the reduced Brillouin zone $\Gamma$MN for two photonic crystals is plotted on the right column. See Section~\ref{s:ex1}.}
% \label{Fig:square_lattice_init_band}
% \end{figure}

% \begin{figure}[!htbp]
% \begin{centering}
% \includegraphics[width=12cm {sqr_joint_PC_shift=0.jpg}
% \includegraphics[width=5.5cm,,trim = 0mm 5mm 0mm 5mm, clip]{figs/square_lattice_Berry_curvature_band1-3_init_V1.jpg}
% \hspace*{-20pt}
% \includegraphics[width=5.5cm,trim = 0mm 5mm 0mm 5mm, clip]{figs/square_lattice_Berry_curvature_band1-3_init_V2.jpg}
% \par\end{centering}
% \caption{The gap Berry curvature $F_{g,3}(\bkappa)$ for two initial photonic crystals shown in Figure \ref{Fig:square_lattice_init_PCs}. The peak values of $F_{g,3}(\bkappa)$ appear at $\bkappa_1=(-1.2, 1.2)$ and $\bkappa_2=(1.2, -1.2)$. See Section~\ref{s:ex1}.}
% \label{Fig:Berry_ini_ex1}
% \end{figure}

% \begin{figure}[!p]
% \begin{centering}
% \includegraphics[scale=0.9,trim = 0mm 5mm 0mm 5mm, clip]{figs/square_final_V.png} 
% \par\end{centering}
% \vspace*{-20pt}
% \caption{Distribution of the electric permittivity for the optimized photonic crystals and the corresponding band structures; see Section~\ref{s:ex1}..}
% \label{Fig:square_lattice_opt_PCs}
% \end{figure}

% \begin{figure}[!htbp]
% \begin{centering}
% \includegraphics[scale=0.65,trim = 0mm 5mm 0mm 5mm, clip]{figs/square_final_bandgap.png} 
% \par\end{centering}
% \caption{The third (left column) and fourth (middle column) band of two photonic crystals in Figure \ref{Fig:square_lattice_opt_PCs}. The eigenvalues along the boundary of the reduced Brillouin zone $\Gamma$MN for two photonic crystals is plotted on the right column. See Section~\ref{s:ex1}.}
% \label{Fig:square_lattice_opt_band}
% \end{figure}

% \begin{figure}[!htbp]
% \begin{centering}
% \includegraphics[width=12cm {sqr_joint_PC_shift=0.jpg}
% \includegraphics[width=5.5cm,trim = 0mm 5mm 0mm 5mm, clip]{figs/square_lattice_Berry_curvature_band1-3_opt_V1.jpg}
% \hspace*{-20pt}
% \includegraphics[width=5.5cm,trim = 0mm 5mm 0mm 5mm, clip]{figs/square_lattice_Berry_curvature_band1-3_opt_V2.jpg}
% \par\end{centering}
% \caption{The gap Berry curvature $F_{g,3}(\bkappa)$ for two optimized photonic crystals.
% The peak values of Berry curvature appear at $\bkappa_1=(-1.2, 1.2)$ and $\bkappa_2=(1.2, -1.2)$; see Section~\ref{s:ex1}. }
% \label{Fig:Berry_opt_ex1}
% \end{figure}

% \begin{figure}[!p]
% \begin{centering}
% \includegraphics[width=16cm]{figs/square_lattice_joint_PC.jpg}
% \includegraphics[width=8cm]{figs/square_lattice_edge_mode_dispersion_opt.jpg}
% \par\end{centering}
% \caption{Top: The joint structure formed by gluing two optimized periodic media together along the $x_2$-axis. Bottom: The dispersion relation of the edge modes along the interface for the wavenumber $\kappa_\parallel\in[-\pi,\pi]$. The gray area is the projection of the continuous spectrum of the periodic media along the interface direction. See Section~\ref{s:ex1}.} 
% \label{Fig:square_lattice_joint_PC}
% \end{figure}

% \begin{figure}[!htbp]
% \begin{centering}
% \includegraphics[width=16cm]{figs/square_lattice_edge_mode_k=0_opt.jpg}
% \includegraphics[width=16cm]{figs/square_lattice_edge_mode1_k=0.6pi_opt.jpg}
% \includegraphics[width=16cm {figs/square_lattice_edge_mode2_k=0.6pi_opt.jpg}
% \par\end{centering}
% \caption{Real part of the edge modes when $\kappa_\parallel=0$ and $0.6\pi$ respectively. The corresponding eigenfrequencies are marked as black dots in the bottom of Figure \ref{Fig:square_lattice_joint_PC}. See Section~\ref{s:ex1}. }
% \label{Fig:square_lattice_edge_mode}
% \end{figure}

\clearpage

\subsection{Hexagonal lattice with C3 symmetry} 
\label{s:ex2}
\begin{figure}[!p]
\begin{centering}
\includegraphics[scale=1,trim = 0mm 5mm 0mm 5mm, clip]{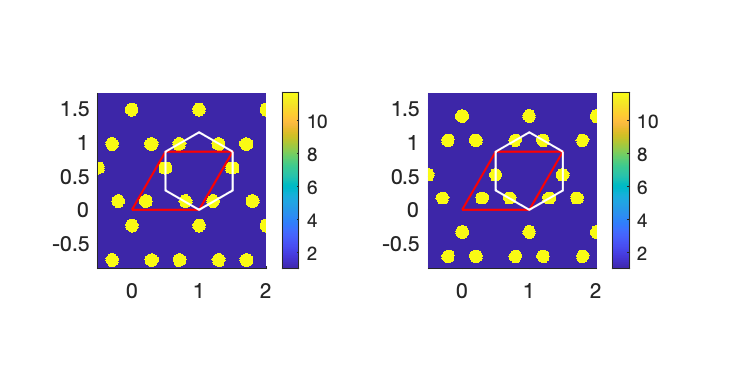}
\includegraphics[scale=0.65]{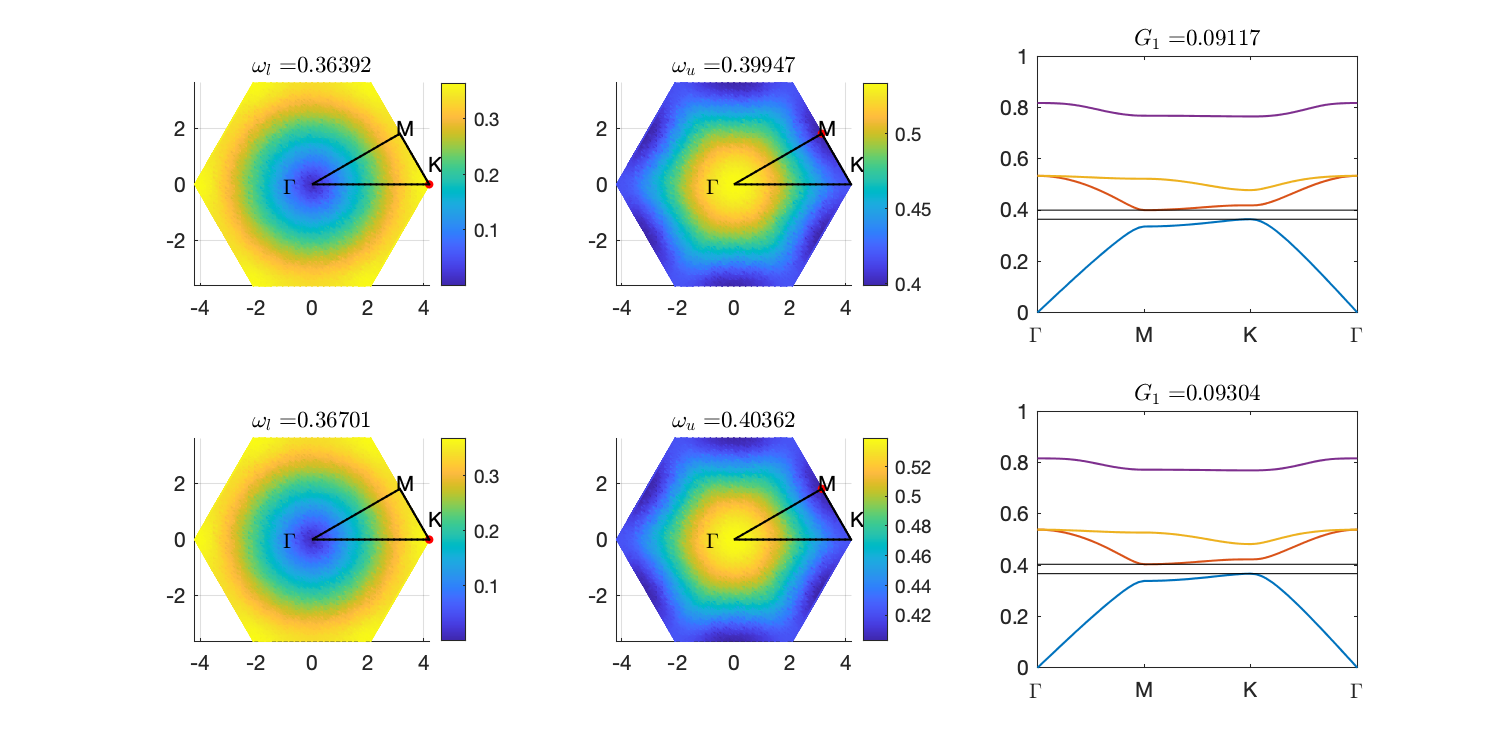}
\includegraphics[width=6cm]{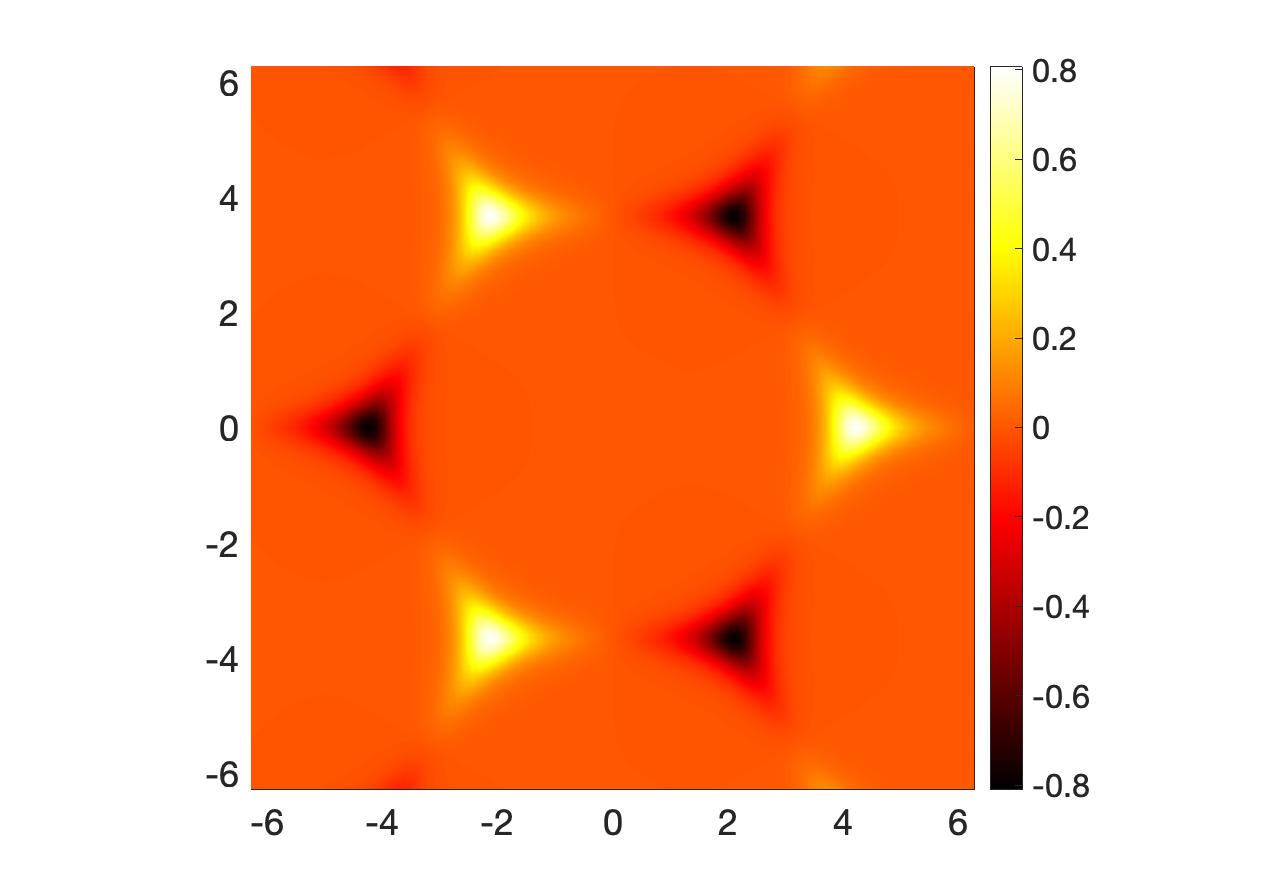}
\includegraphics[width=6cm]{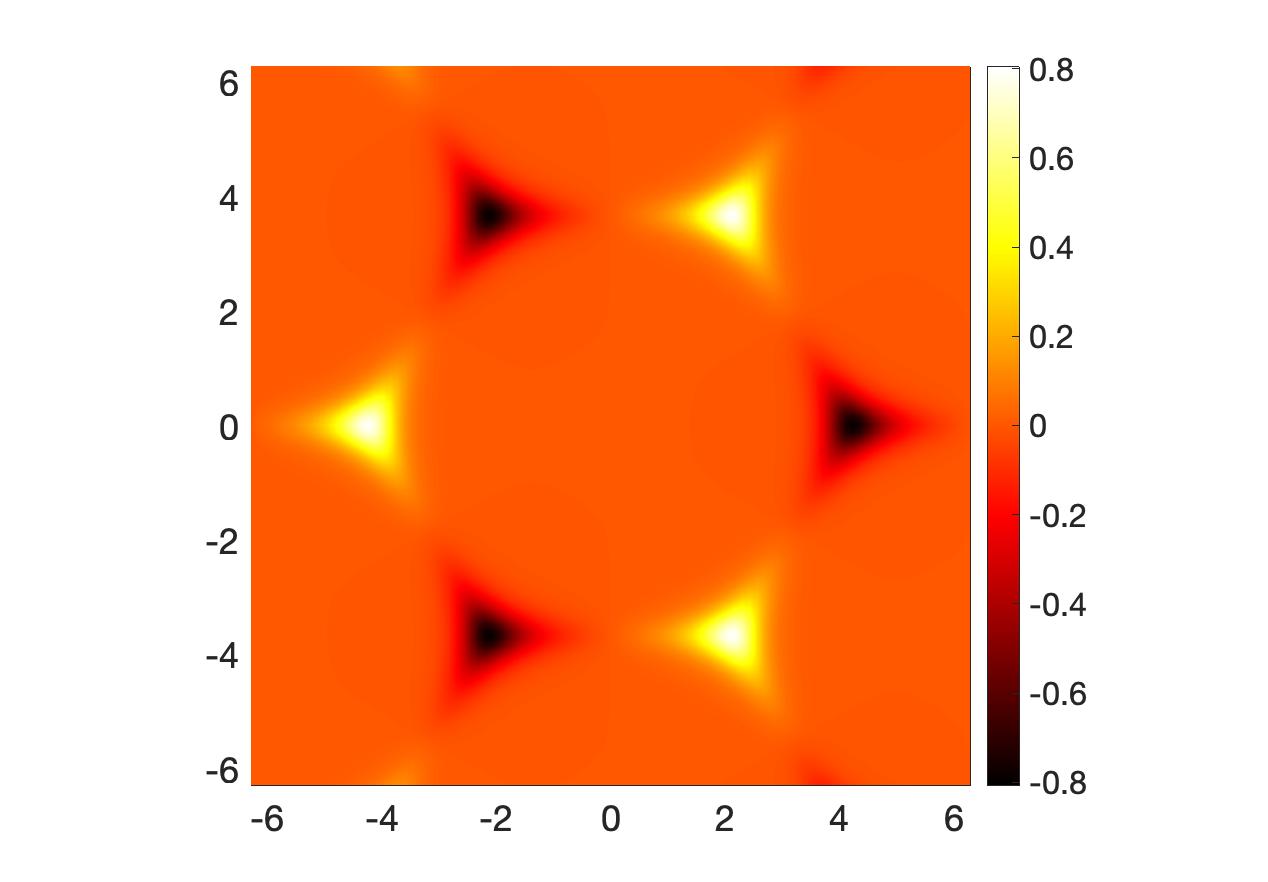}
\end{centering}
\caption{Row 1: Distribution of the electric permittivity for the initial configuration of photonic crystals. The three cylindrical rods in the fundamental periodic cell are pull away (left) or closer (right) to each other along the hexagonal edges by a distance of $0.05$. The parallelogram represents the fundamental periodic cell of the photonic crystal.
Row 2-3: The first (left column) and second (middle column) band of the photonic crystals. The eigenvalues of two PCs along the boundary of the reduced Brillouin zone $\Gamma$MK is plotted in the right column.
Row 4: The gap Berry curvature $F_{g,1}(\bkappa)$ for the initial photonic crystals. The peak values appear at the corners of the Brillouin zone.
See Section~\ref{s:ex2}. }
\label{Fig:tri_C3_init_PCs}
\end{figure}

In this example, we consider PCs with the periodicity of a hexagonal lattice, for which the lattice vectors are given by $\be_1=(1,0)$ and $\be_2=\left( \frac{1}{2},\frac{\sqrt{3}}{2}\right)$. 
Figure \ref{Fig:tri_C3_init_PCs} shows two initial PCs consist of cylindrical rods of diameter  $0.2$. 
The positions of the three rods in the fundamental periodic cells (the parallelogram in Figure \ref{Fig:tri_C3_init_PCs}) differ between the two PCs; in one PC the rods are pulled away from each other along the hexagonal edges by a distance of $0.05$. This opens a small gap near a Dirac point, which is formed by the first two bands of the spectrum when the rods are located in the middle of the hexagonal edges \cite{wong2020gapless}. 

The initial band structure is shown in Figure \ref{Fig:tri_C3_init_PCs} (Row 2-3), and the initial shared gap-to-midgap ratio between the first two bands is $G_1(\varepsilon_1,\varepsilon_2)\approx 0.0827$ and the individual ratios $G_1(\varepsilon_1)$ and $G_1(\varepsilon_2)$ are indicated in Figure \ref{Fig:tri_C3_init_PCs}.    
The corresponding gap Berry curvature $F_{g,1}(\bkappa)$ for two PCs is plotted in Figure \ref{Fig:tri_C3_init_PCs} (Row 4), where the peak values appear at the corners of the hexagonal Brillouin zone. 
By choosing $\bkappa_1=(\frac{4\pi}{3}, 0)$ and $\bkappa_2=\frac{2\pi}{3}(1,\sqrt{3})$, the valley Chern number takes the value $1$ and $-1$ for the two photonic crystals respectively.

\begin{figure}[!p]
\begin{centering}
\par\end{centering}
\includegraphics[scale=1,trim = 0mm 5mm 0mm 5mm, clip]{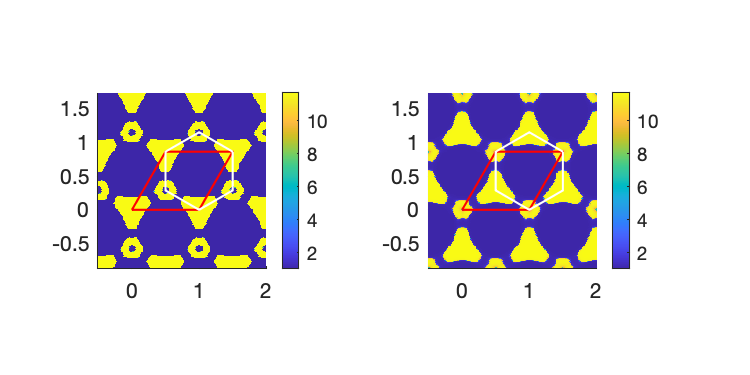}
\includegraphics[scale=0.65]{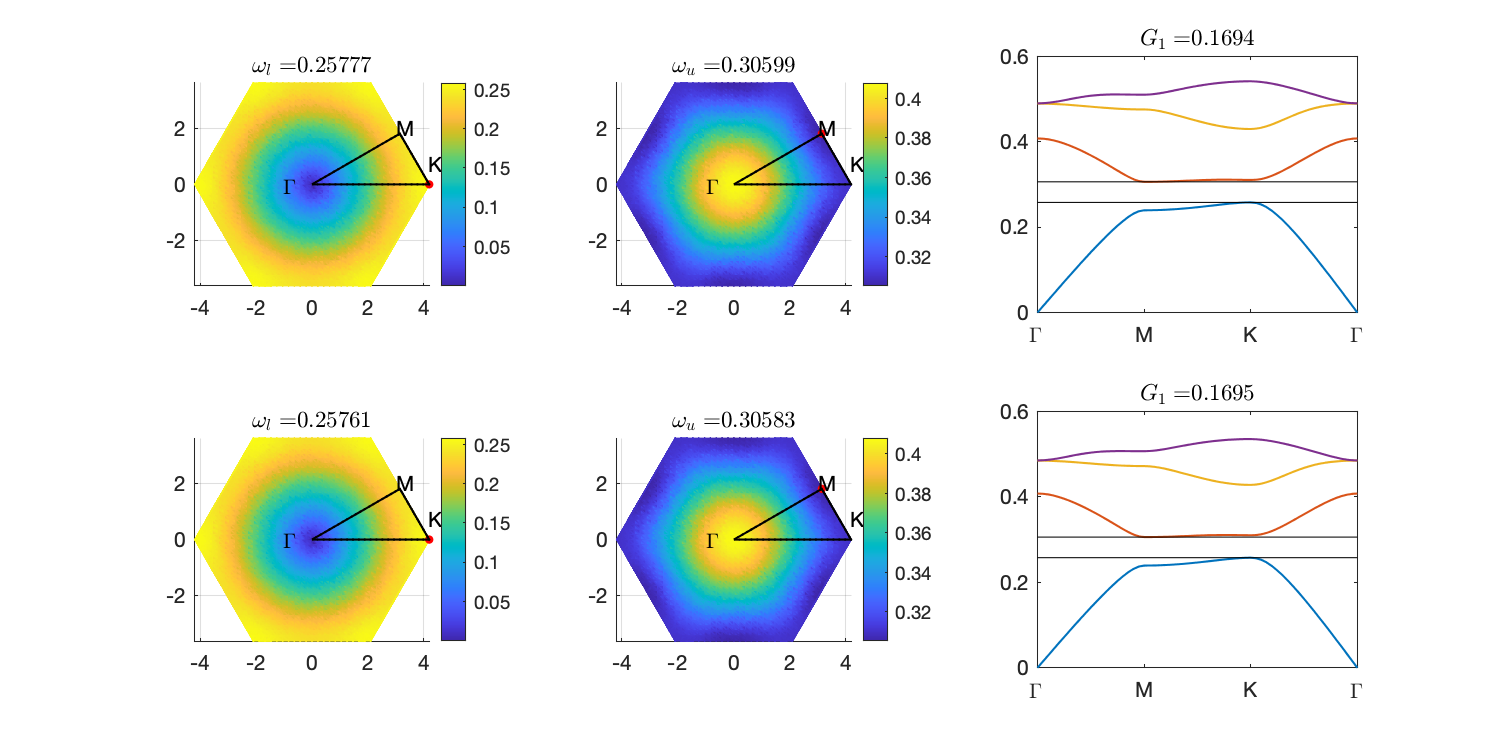}
\includegraphics[width=6.5cm]{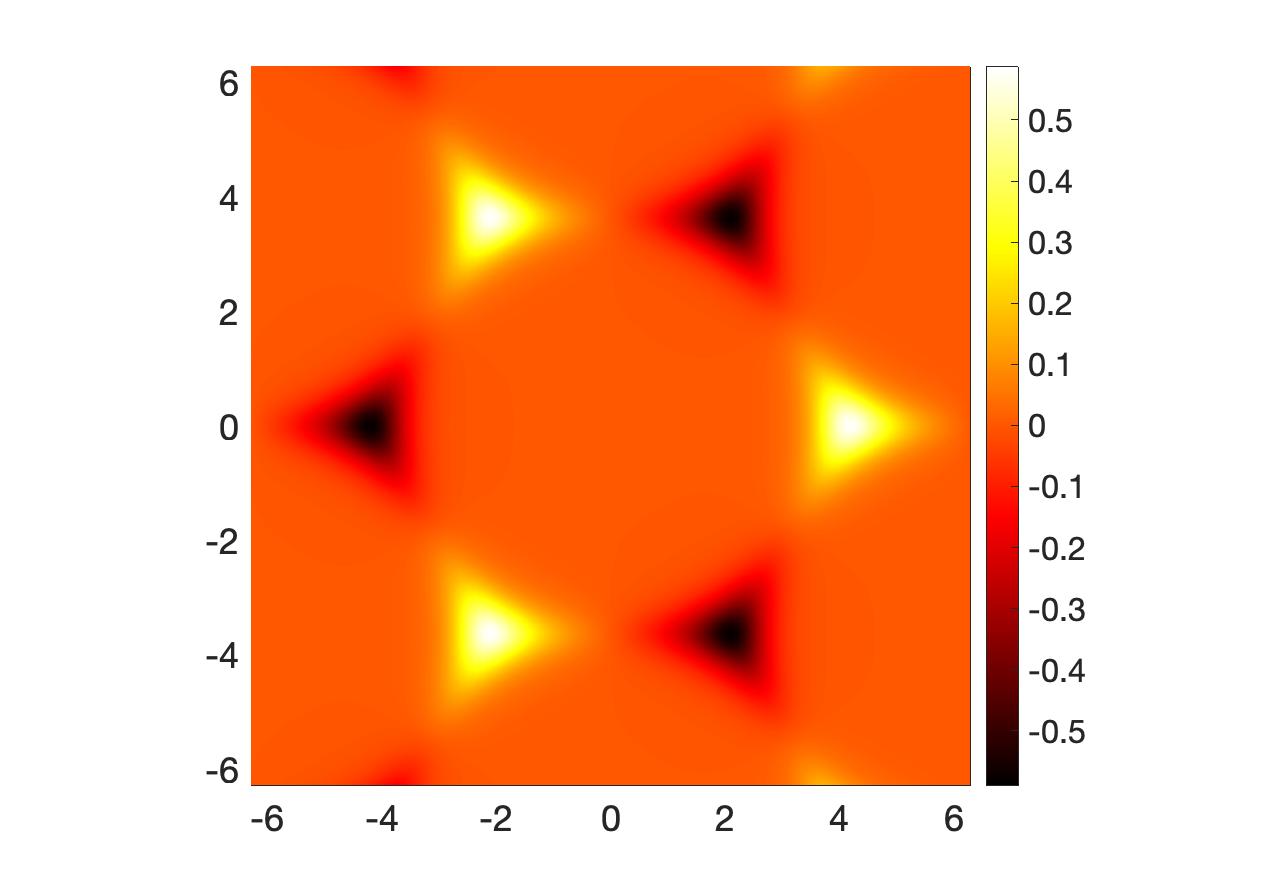}
\hspace*{-20pt}
\includegraphics[width=6.5cm]{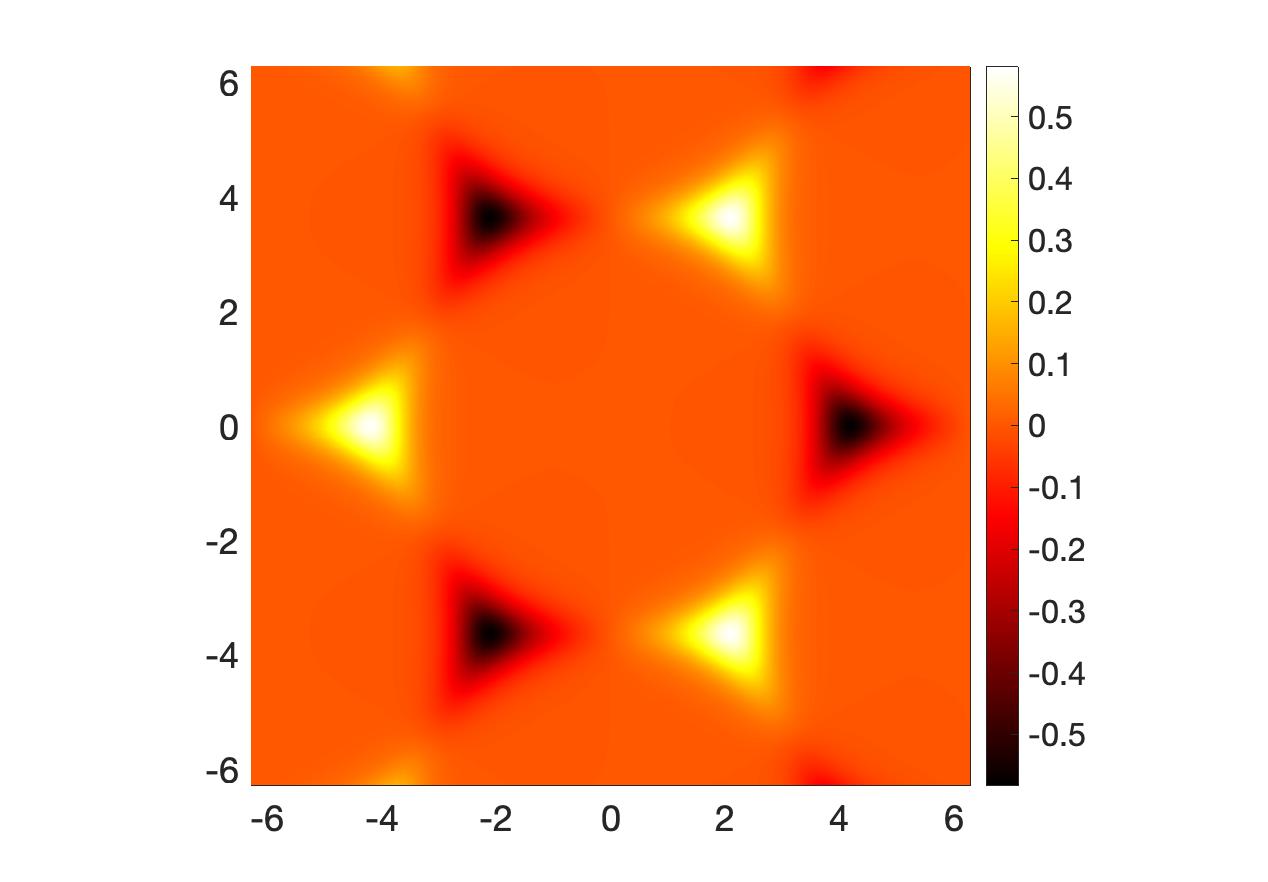}
\caption{Row 1: Distribution of the electric permittivity for the optimized photonic crystals.
Row 2-3: The first (left column) and second (middle column) band of the optimized photonic crystals in Figure \ref{Fig:tri_C3_opt_PCs}. The eigenvalues of two PCs along the boundary of the reduced Brillouin zone $\Gamma$MK is plotted in the right column.
Row 4: The gap Berry curvature $F_{g,1}(\bkappa)$ for the optimized photonic crystals.
See Section~\ref{s:ex2}.}
\label{Fig:tri_C3_opt_PCs}
\end{figure}

We preserve the $C_3$ symmetry of the hexagonal lattice in the optimization. Furthermore, similar to the square lattice, we set $\tau_0=\frac{1}{2}$ in \eqref{eq:constraint_F} to preserve the valley Chern number of the periodic media so that the Berry curvature does not diminish during the iteration. The optimized PCs are obtained after 5 iterations and plotted in Figure~\ref{Fig:tri_C3_opt_PCs}. 
The first and second dispersion surfaces are plotted in Figure~\ref{Fig:tri_C3_opt_PCs} (Row 2-3) and the shared gap-to-midgap ratio is $G_1\approx 0.1689$. 
The gap Berry curvature $F_{g,1}(\bkappa)$ for the optimized PCs is plotted in Figure~\ref{Fig:tri_C3_opt_PCs} (Row 4), which still attains the peak values at the corners of the Brillouin zone as the initial structures and the valley Chern number is preserved.

\begin{figure}[!p]
\begin{centering}
\hspace*{-30pt}
\includegraphics[width=18cm]{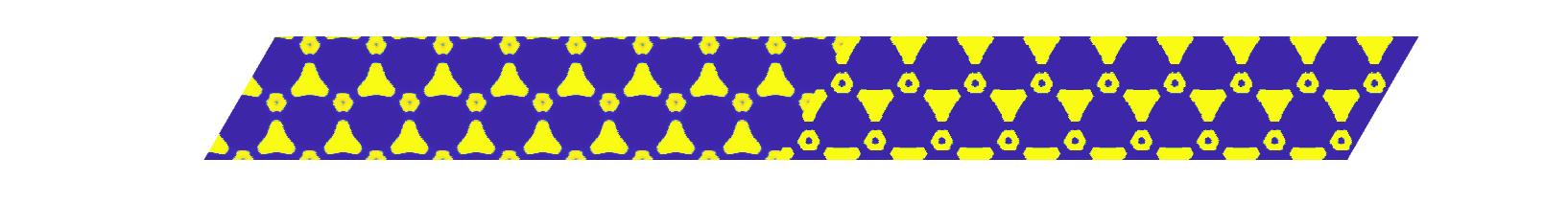}
\includegraphics[width=8cm]{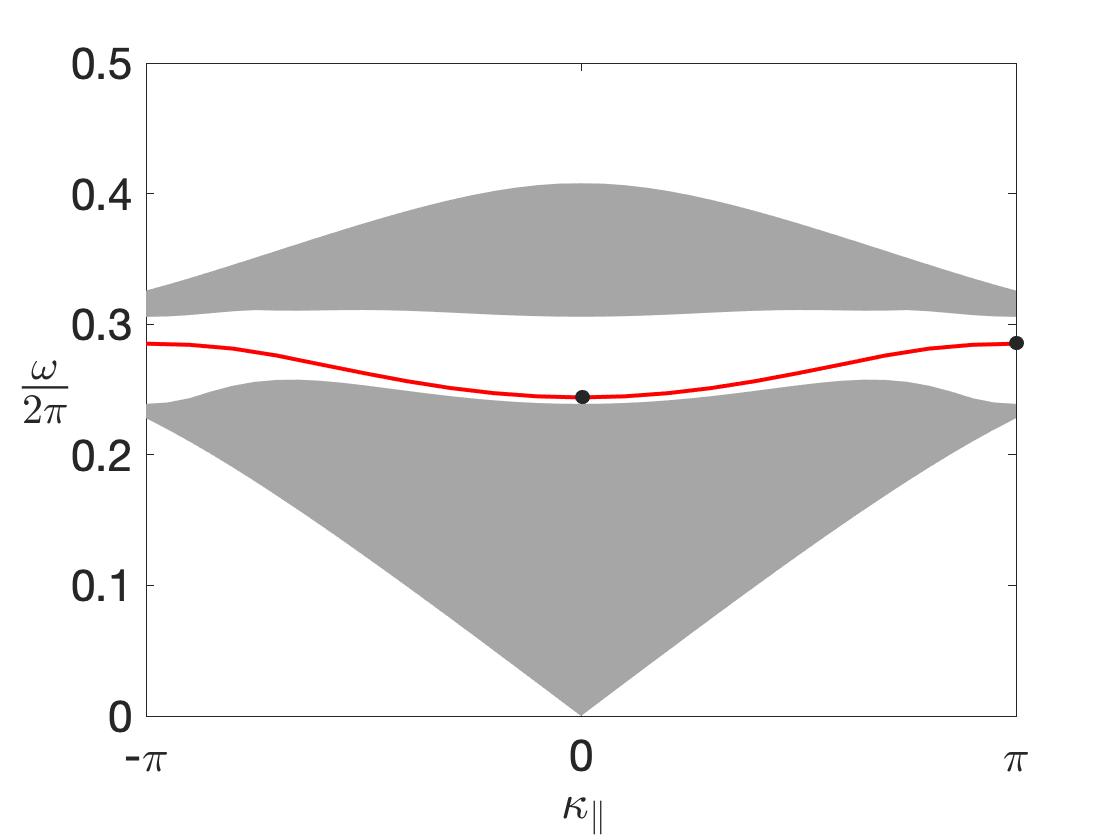}
\hspace*{-30pt}
\includegraphics[width=18cm]{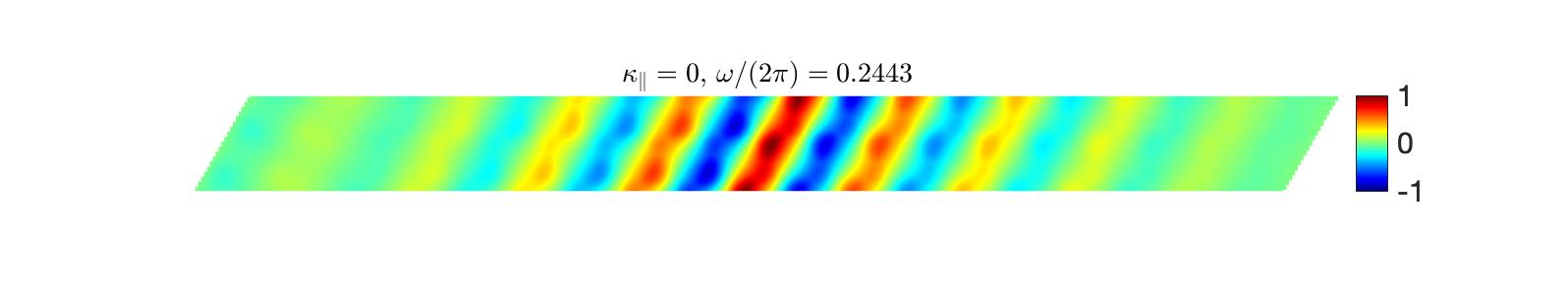}\\
\vspace*{-30pt}
\hspace*{-30pt}
\includegraphics[width=18cm]{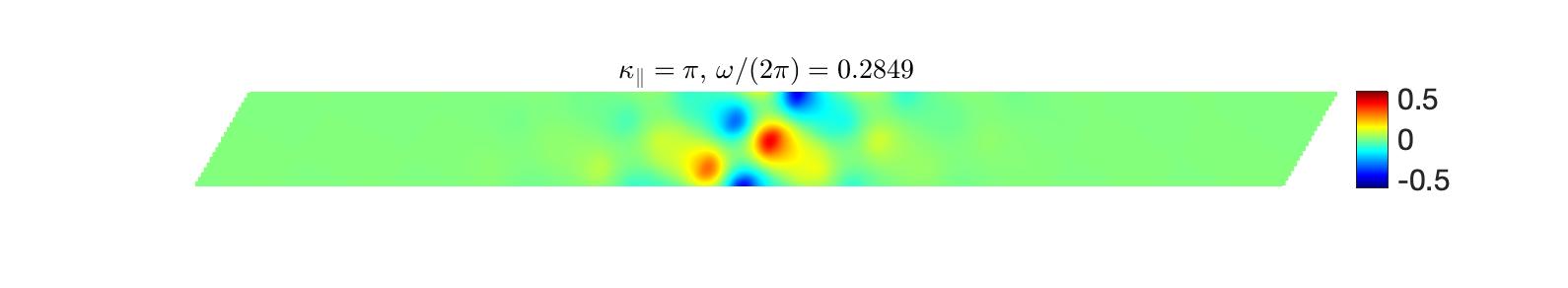} \\
\vspace*{-20pt}
\par\end{centering}
\caption{Row 1: The joint structure formed by gluing two optimized periodic media together along the $\be_2$ direction. Row 2: The dispersion relation of the edge modes along the interface for the wavenumber $\kappa_\parallel\in[-\pi,\pi]$. The gray area is the projection of the continuous spectrum of the periodic media along the $\be_2$ direction.
Row 3-4: Real part of the edge modes when $\kappa_\parallel=0$ and $\pi$ respectively. The corresponding eigenfrequencies are marked as black dots in the dispersion relation.
See Section~\ref{s:ex2}.}
\label{Fig:tri_C3_edge_mode_dispersion_opt}
\end{figure}

To compute the bandwidth of the edge modes, we glue the two optimized photonic crystals together along the interface direction $\be_2$ as shown in Figure~\ref{Fig:tri_C3_edge_mode_dispersion_opt} (Row 1). 
The dispersion relation of the edge modes for $\kappa_\parallel\in[-\pi, \pi]$ is plotted in Figure~\ref{Fig:tri_C3_edge_mode_dispersion_opt} (Row 2), which demonstrates an enlarged bandwidth for the edge modes between the first and the second band of the periodic media. We also plot the edge modes, which is localized along interface of the joint structure,
when $\kappa_\parallel=0$ and $\pi$ in Figure~\ref{Fig:tri_C3_edge_mode_dispersion_opt} (Row 3-4).

% \begin{figure}[!p]
% \begin{centering}
% \hspace*{-30pt}
% \includegraphics[width=18cm]{figs/tri_lattice_C3_joint_PC.jpg}
% \includegraphics[width=8cm]{figs/tri_lattice_C3_dispersion_opt_structure.jpg}
% \par\end{centering}
% \caption{Top: The joint structure formed by gluing two optimized periodic media together along the $\be_2$ direction; Bottom: The dispersion relation of the edge modes along the interface for the wavenumber $\kappa_\parallel\in[-\pi,\pi]$. The gray area is the projection of the continuous spectrum of the periodic media along the $\be_2$ direction. See Section~\ref{s:ex2}.}
% \label{Fig:tri_C3_edge_mode_dispersion_opt}
% \end{figure}

% \begin{figure}[!htbp]
% \begin{centering}
% \hspace*{-30pt}
% \includegraphics[width=18cm]{figs/tri_lattice_C3_edge_mode_opt_k=0.jpg}
% \hspace*{-30pt}
% \includegraphics[width=18cm]{figs/tri_lattice_C3_edge_mode_opt_k=pi.jpg}
% \par\end{centering}
% \caption{Real part of the edge modes when $\kappa_\parallel=0$ and $\pi$ respectively. The corresponding eigenfrequencies are marked as black dots in the bottom of Figure \ref{Fig:tri_C3_edge_mode_dispersion_opt}. See Section~\ref{s:ex2}.}
% \label{Fig:tri_C3_edge_mode}
% \end{figure}

\clearpage
\subsection{Hexagonal lattice with C6 symmetry}
\label{s:ex3}
\begin{figure}[!p]
\begin{centering}
\includegraphics[scale=1,trim = 0mm 5mm 0mm 5mm, clip]{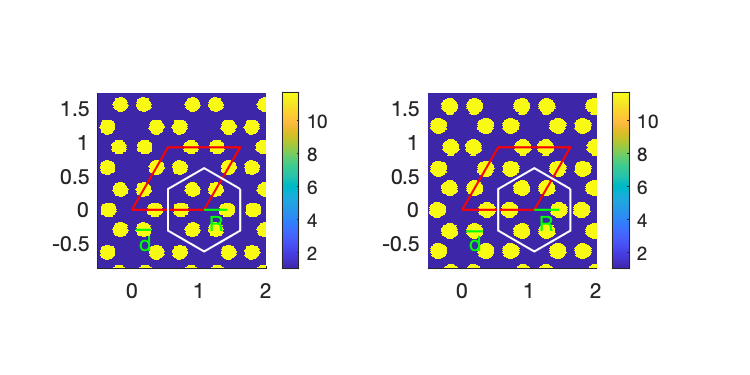}
\includegraphics[scale=0.35,trim = 0mm 5mm 0mm 5mm, clip]{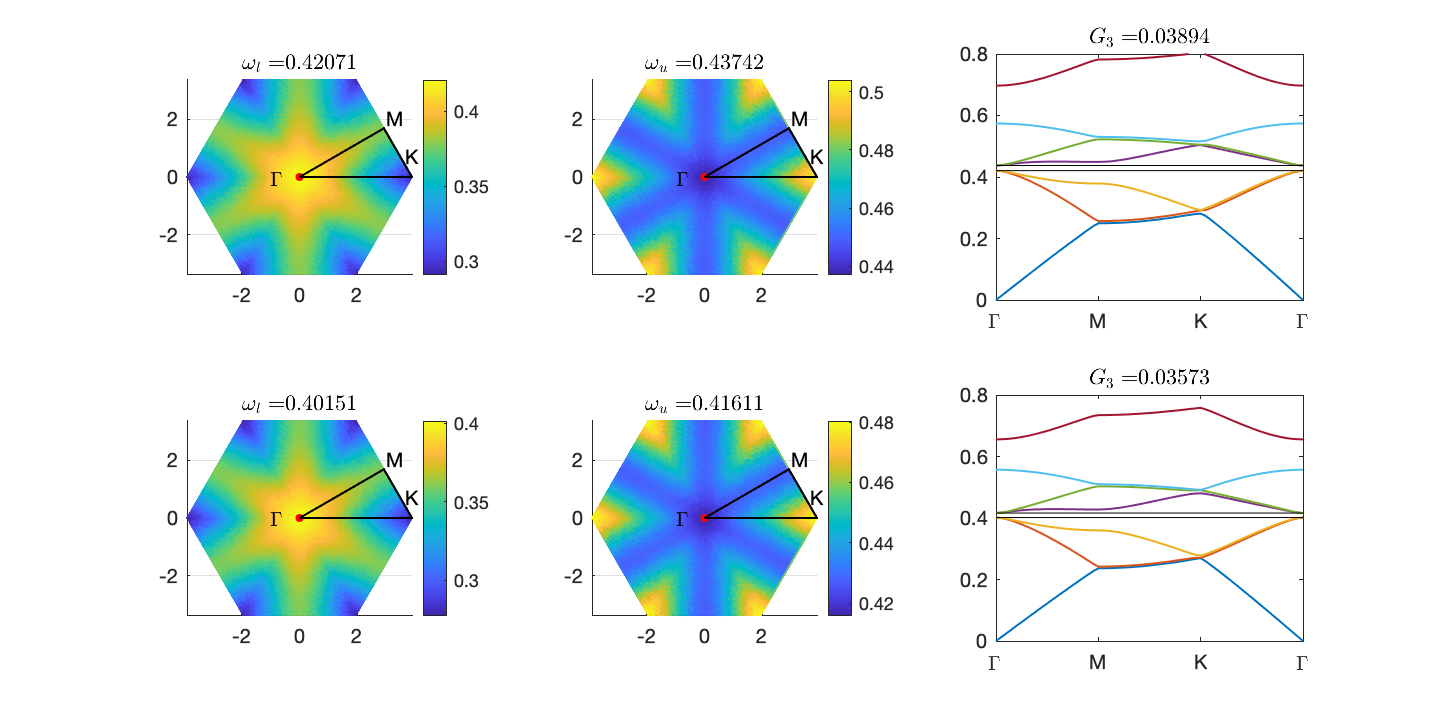}
\includegraphics[width=6cm,trim = 0mm 5mm 0mm 5mm, clip]{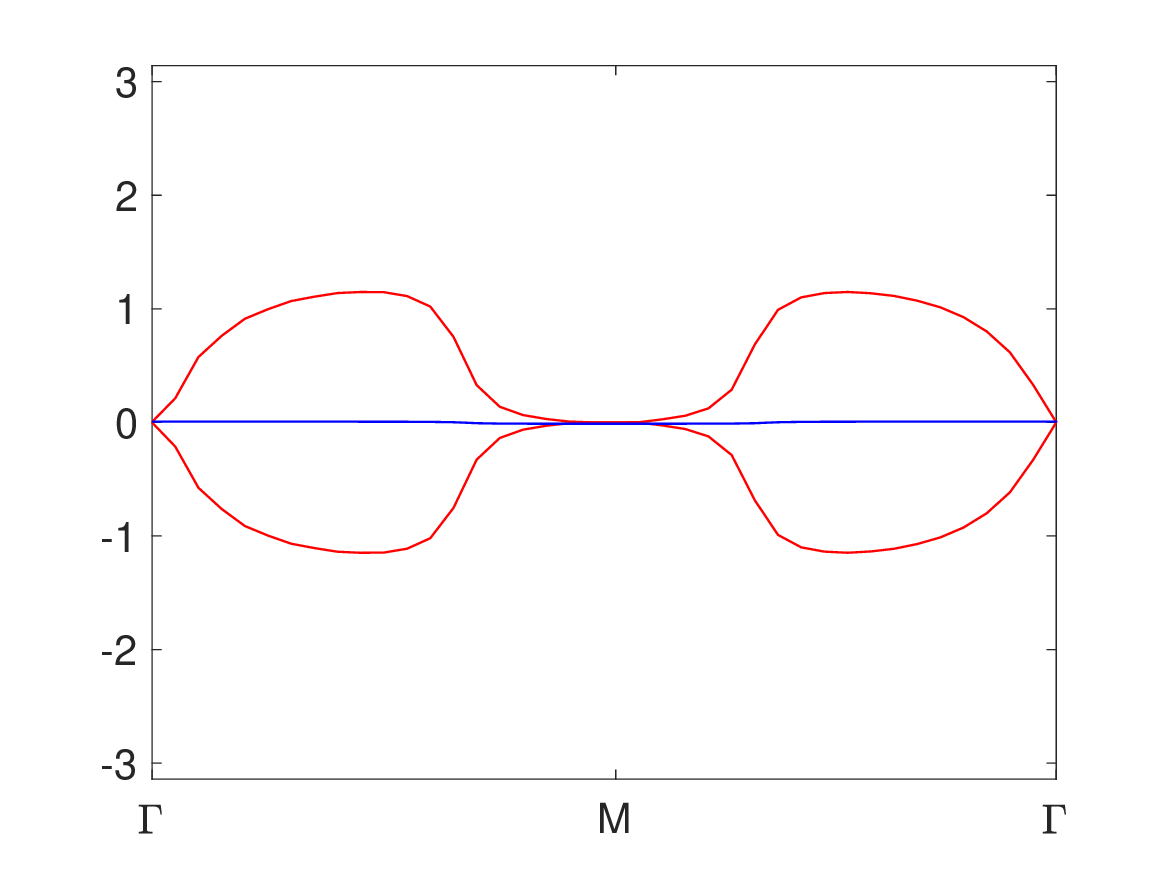}
\includegraphics[width=6cm,trim = 0mm 5mm 0mm 5mm, clip]{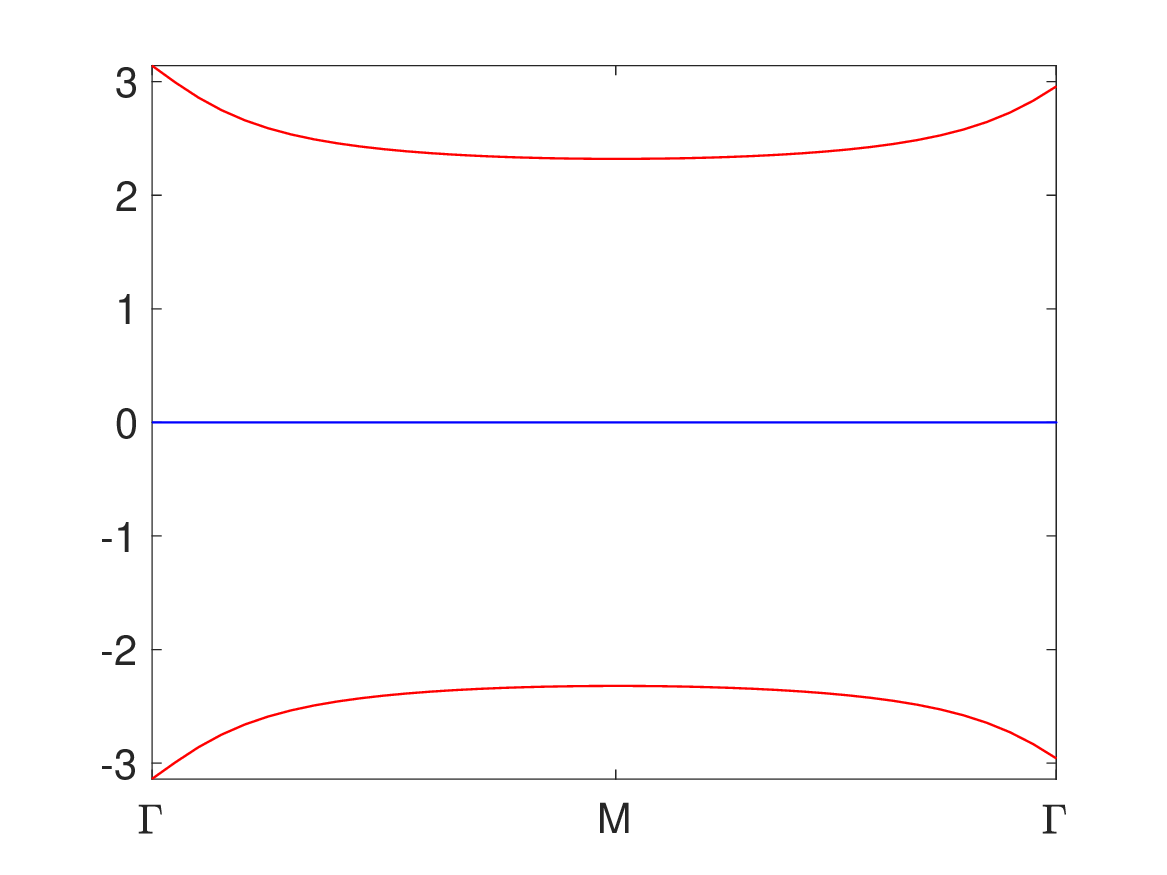}
\par\end{centering}
\caption{Row 1: Distribution of the electric permittivity for the initial configuration of photonic crystals.
Row 2-3: The third (left column) and fourth band (middle column) of two photonic crystals. The eigenvalues along the boundary of the reduced Brillouin zone $\Gamma$MK for two photonic crystals is plotted on the right column. 
Row 4: Wilson loops for the first three bands of the two photonic crystals. Left: the Wilson loops for the first photonic crystal. The value for the first band (blue) is computed via \eqref{eq:wilson_lopp1}, and the value for the degenerate second and third band (red) is computed via \eqref{eq:wilson_lopp2} with $m=2$. Right: the Wilson loops for the second photonic crystal, which is computed via \eqref{eq:wilson_lopp2} with $m=3$ for the first three degenerate bands. 
See Section~\ref{s:ex3}.}
\label{Fig:tri_C6_init_PCs}
\end{figure}

In this example, we consider PCs with the periodicity of a hexagonal lattice and $C_6$ symmetry. Figure \ref{Fig:tri_C6_init_PCs} shows two initial PCs where the fundamental hexagonal periodic cell intersects six cylindrical rods   \cite{wu2015scheme}. Let $d$ be the diameter of cylindrical rod and $R$ be distance between the hexagonal center and the axis of one of the cylindrical rods. Then the rods in the two photonic crystals attain the diameter $d=\frac{2}{3}R$, with the distance given by $R = 1/3.125$ and $R = 1/2.9$ respectively. 

The tiny spectral gap is opened near the double Dirac cones at $\Gamma$ point formed by the second to the fifth band when the distance $R$ between the hexagonal center and the axis of the rod is $1/3$ \cite{wu2015scheme}. The corresponding initial band structure is shown in Figure \ref{Fig:tri_C6_init_PCs} (Row 2-3).  

Here the photonic crystals mimic the spin Hall effect of electron models. In particular, the pseudo spin-up and spin-down states are represented by anti-clockwise and clockwise circular polarization of in-plane magnetic field $\mathbf{H}$. By using the so-called $\bkappa\cdot \mathbf{P}$ perturbation theory, it can be shown that the effective Hamiltonian shares a similar form as Bernevig-Hughes-Zhang (BHZ) model for quantum well, and the corresponding spin gap Chern number is given by $0$ and $2$ respectively for the two effective Hamiltonians \cite{wu2015scheme}. Alternatively, the topological phase difference of two periodic media can be manifested through the Wilson loop as plotted in Figure \ref{Fig:tri_C6_init_PCs} (Row 4), which is the Berry phase computed in a closed loop along the $\bkappa_2$ direction in the $\bkappa$-space for fixed $\bkappa_1$ (see Appendix). An important property of the Wilson loop is that its spectrum is equivalent to the centers of the localized Wannier functions associated with the group of bands up to a scaling constant \cite{vanderbilt2018berry}. For trivial bands, the Wannier functions are maximally localized, while for systems that possess non-trivial topology, the Wannier functions are delocalized. 
As observed from Figure \ref{Fig:tri_C6_init_PCs}, for the first photonic crystal, the Wilson loop of the first band has a constant value equal to 0, which is representative of a trivial gap. For the set of the second and third bands, the eigenvalues of the Wilson loop associated with the second and third band are localized near 0, showing the trivial character of the periodic medium. For the second photonic crystal, other than the trivial Wilson loop value equal to 0, the other two eigenvalues of Wilson loop are localized near $\pm \pi$, which demonstrates the non-trivial topology of the periodic medium.

\begin{figure}[!p]
\begin{centering}
\includegraphics[scale=1,trim = 0mm 5mm 0mm 5mm, clip]{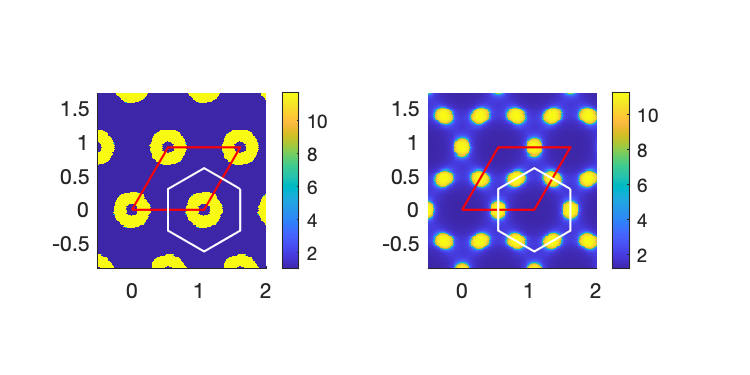}
\includegraphics[scale=0.65,trim = 0mm 5mm 0mm 5mm, clip]{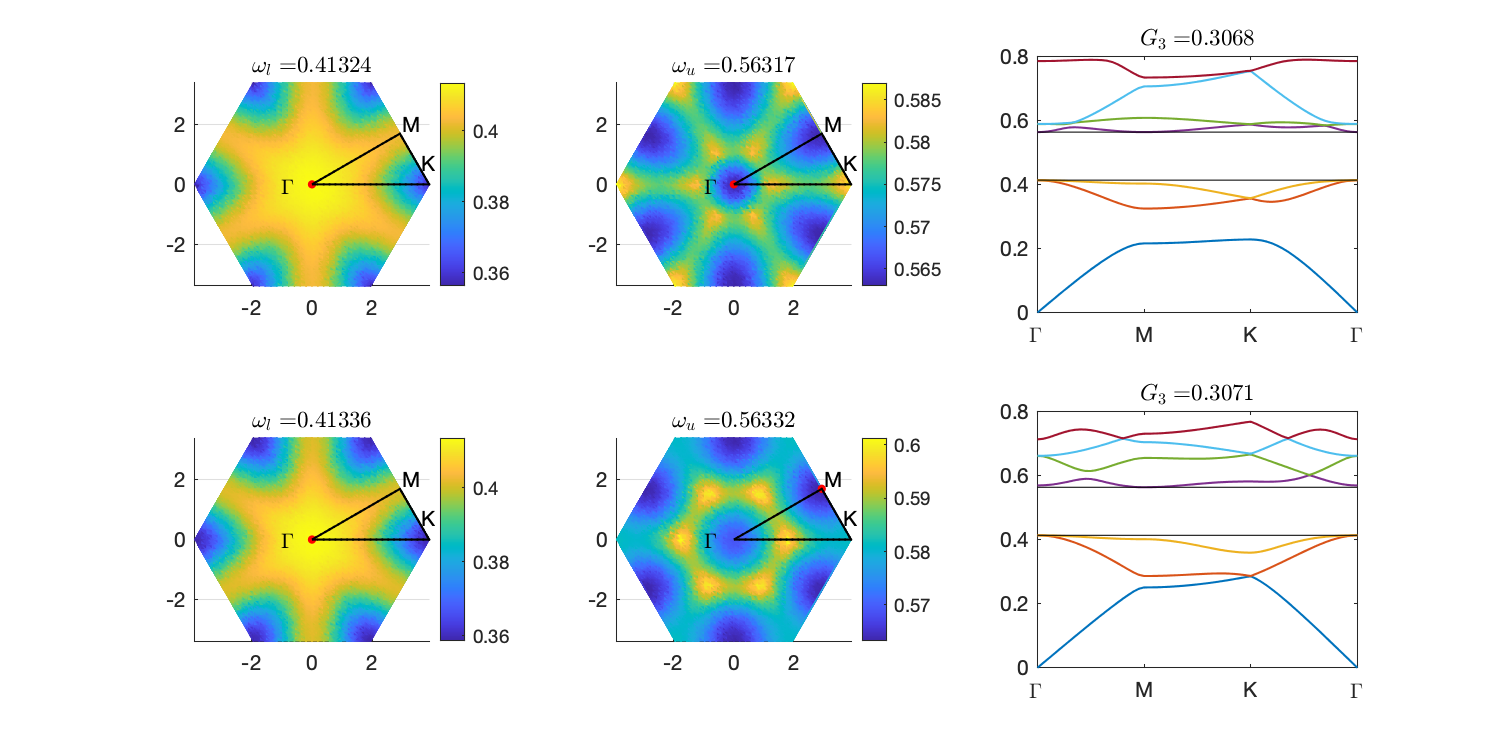}
\includegraphics[width=6cm,trim = 0mm 5mm 0mm 5mm, clip]{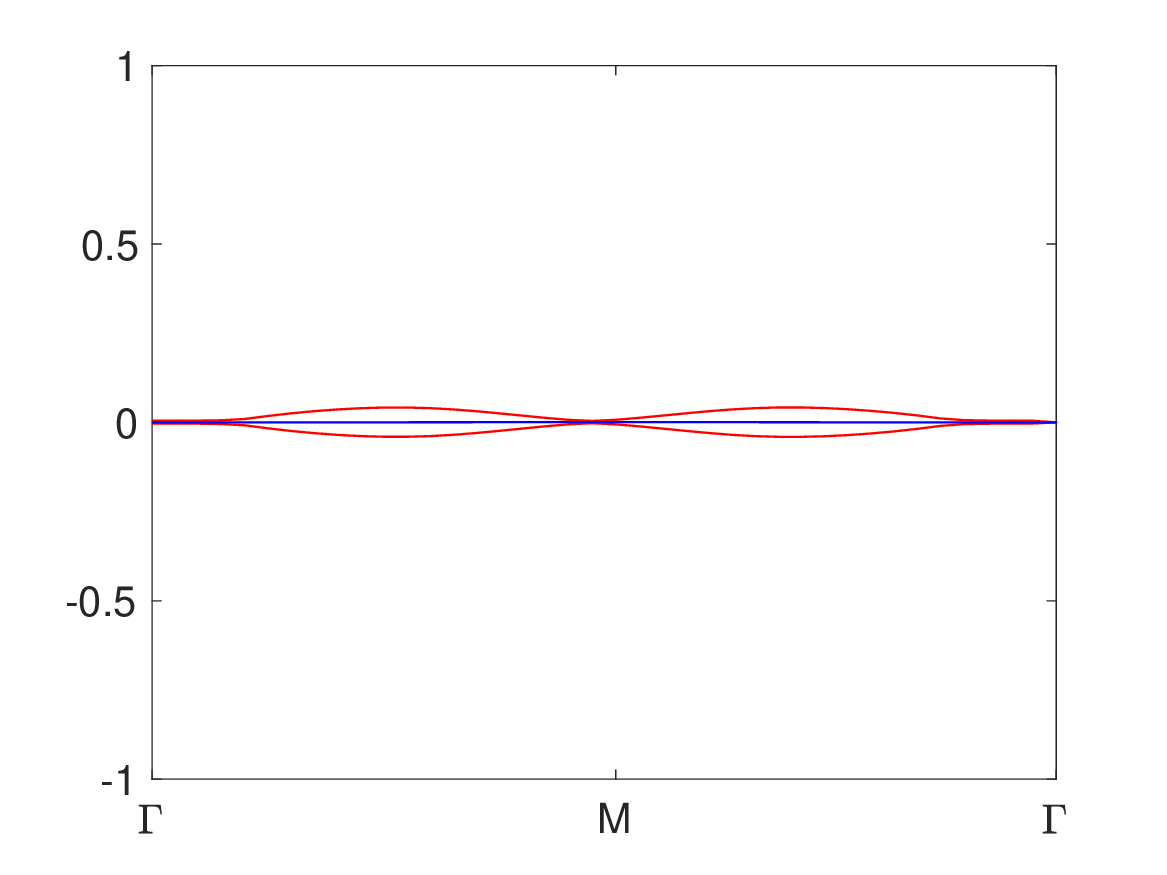}
\includegraphics[width=6cm,trim = 0mm 5mm 0mm 5mm, clip]{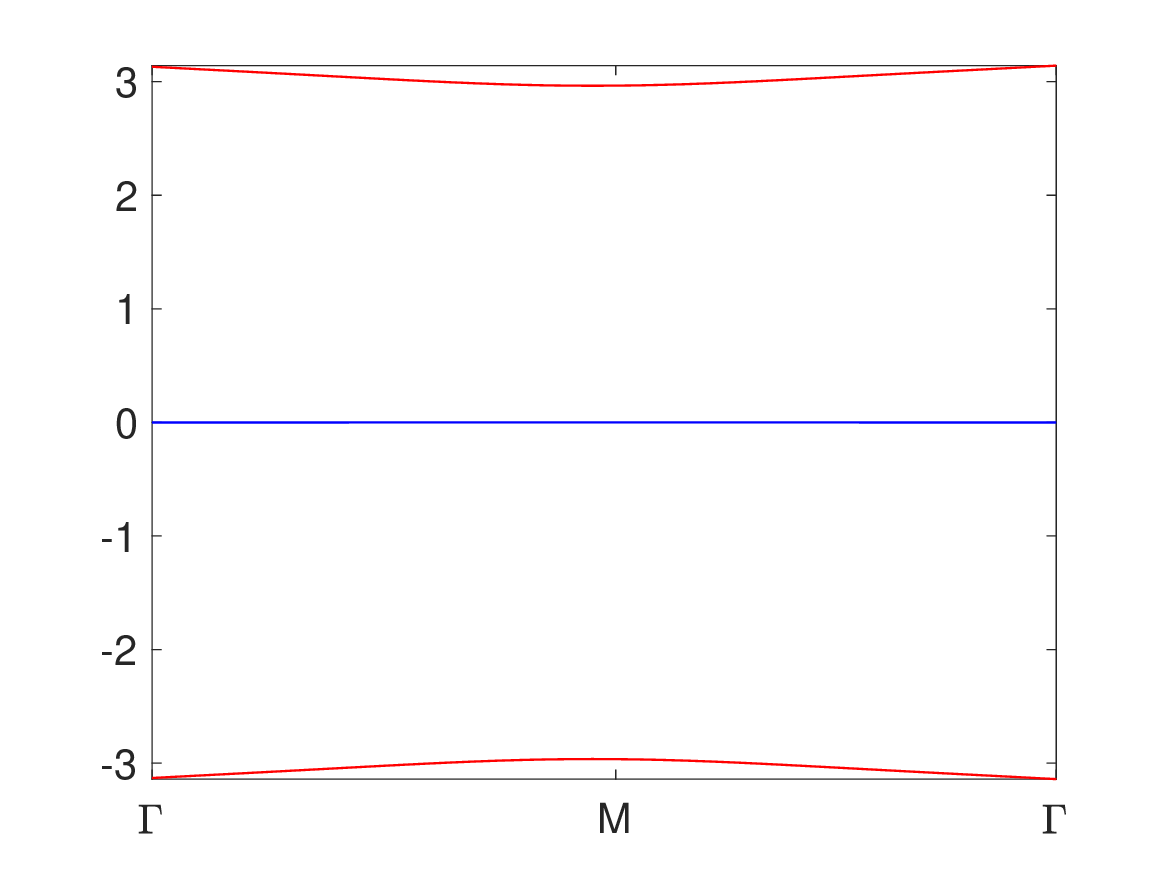}
\par\end{centering}
\caption{Row 1: Distribution of the electric permittivity for the optimized photonic crystals. Row 2-3: The third (left column) and fourth band (middle column) of two photonic crystals. The eigenvalues along the boundary of the reduced Brillouin zone $\Gamma$MK for two photonic crystals is plotted on the right column. Row 4: Wilson loops for the first three bands of the two optimized photonic crystals in Figure \ref{Fig:tri_C6_opt_PCs}. The calculation for each band is the same as in Figure \ref{Fig:tri_C6_init_PCs}.
See Section~\ref{s:ex3}.}
\label{Fig:tri_C6_opt_PCs}
\end{figure}

The C6 symmetry of the lattice is preserved during the optimization iterations, and we obtain two optimized periodic media shown in Figure \ref{Fig:tri_C6_opt_PCs} after 10 iterations. The corresponding band structure is shown in Figure \ref{Fig:tri_C6_opt_PCs} (Row 2-3), and the shared gap-to-midgap $G_3\approx 0.3066$, which is much larger than the initial shared gap-to-midgap. Since the integer spin Chern number and the Berry phase remain  the same when the medium coefficient evolves unless the spectral gap closes, there is no need to impose the constraint on the topological indices explicitly in optimization framework \eqref{e:OptProb-TM}. As demonstrated in Figure \ref{Fig:tri_C6_opt_PCs} (Row 4), the localization of the Wilson loop for the optimized periodic media remains the same compared to the initial photonic crystals, as such the topological nature of the media does not change after the optimization.

\begin{figure}[!p]
\begin{centering}
\hspace*{-30pt}
\includegraphics[width=18cm]{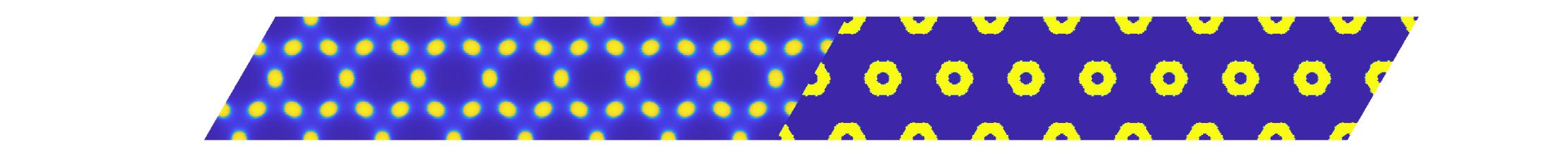}
\includegraphics[width=8cm]{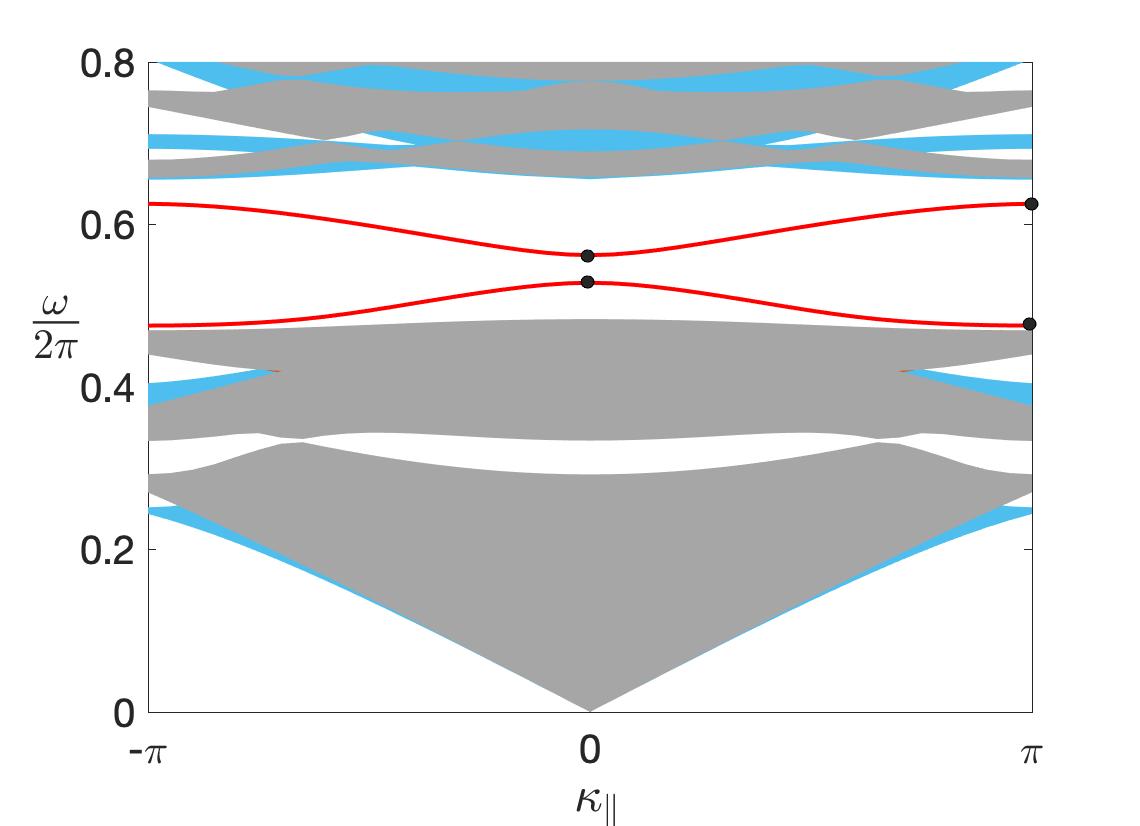}
\hspace*{-30pt}
\vspace*{-30pt}
\includegraphics[width=18cm]{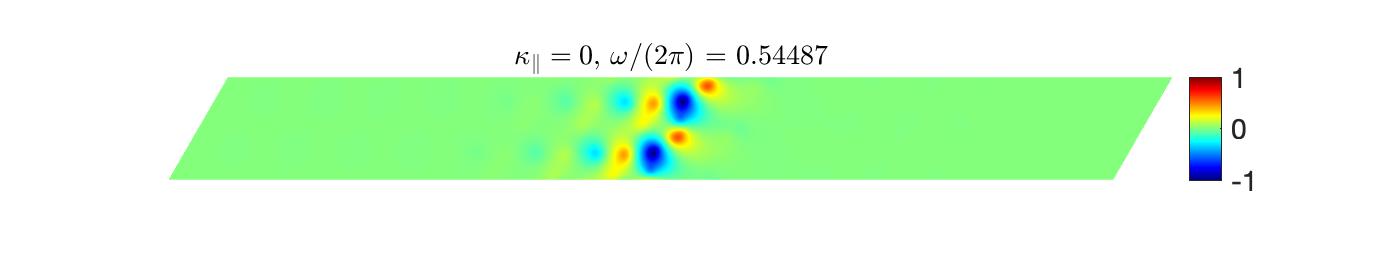}
\hspace*{-30pt}
\vspace*{-30pt}
\includegraphics[width=18cm]{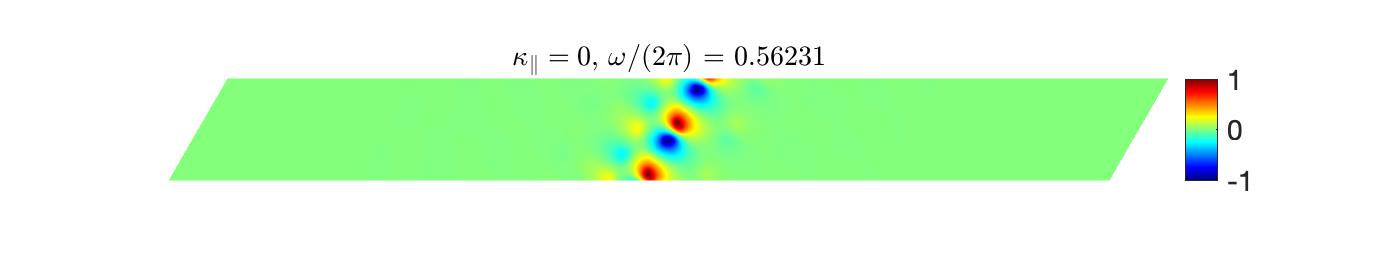}
\hspace*{-30pt}
\vspace*{-30pt}
\includegraphics[width=18cm]{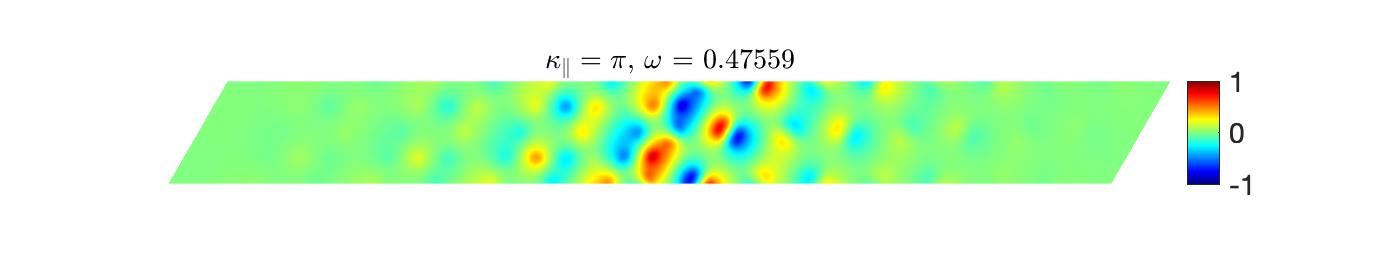}
\hspace*{-30pt}
\vspace*{-30pt}
\includegraphics[width=18cm]{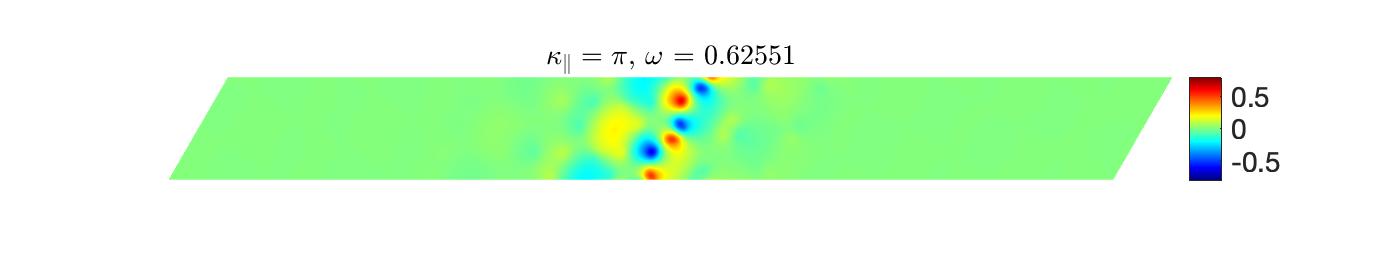}
\par\end{centering}
\caption{Row 1: The joint structure formed by gluing two optimized periodic media together along the $\be_2$ direction. Row 2: The dispersion relation of the edge modes along the interface for the wavenumber $\kappa_\parallel\in[-\pi,\pi]$.
The gray and blue areas denote the projection of the continuous spectrum for two periodic media along the $\be_2$ direction.
Row 3-6: Real part of the edge modes when $\kappa_\parallel=0$ and $\pi$ respectively. The corresponding eigenfrequencies are marked as black dots in the bottom of the dispersion relation.
See Section~\ref{s:ex3}.
}
\label{Fig:tri_C6_edge_mode_dispersion1}
\end{figure}

\begin{figure}[!htbp]
\begin{centering}
\hspace*{-30pt}
\includegraphics[width=18cm]{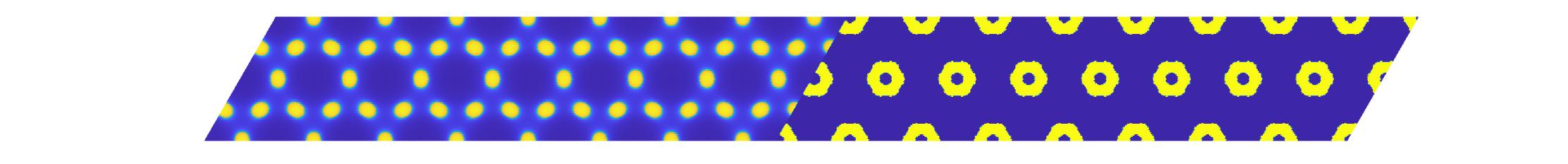}
\includegraphics[width=8cm]{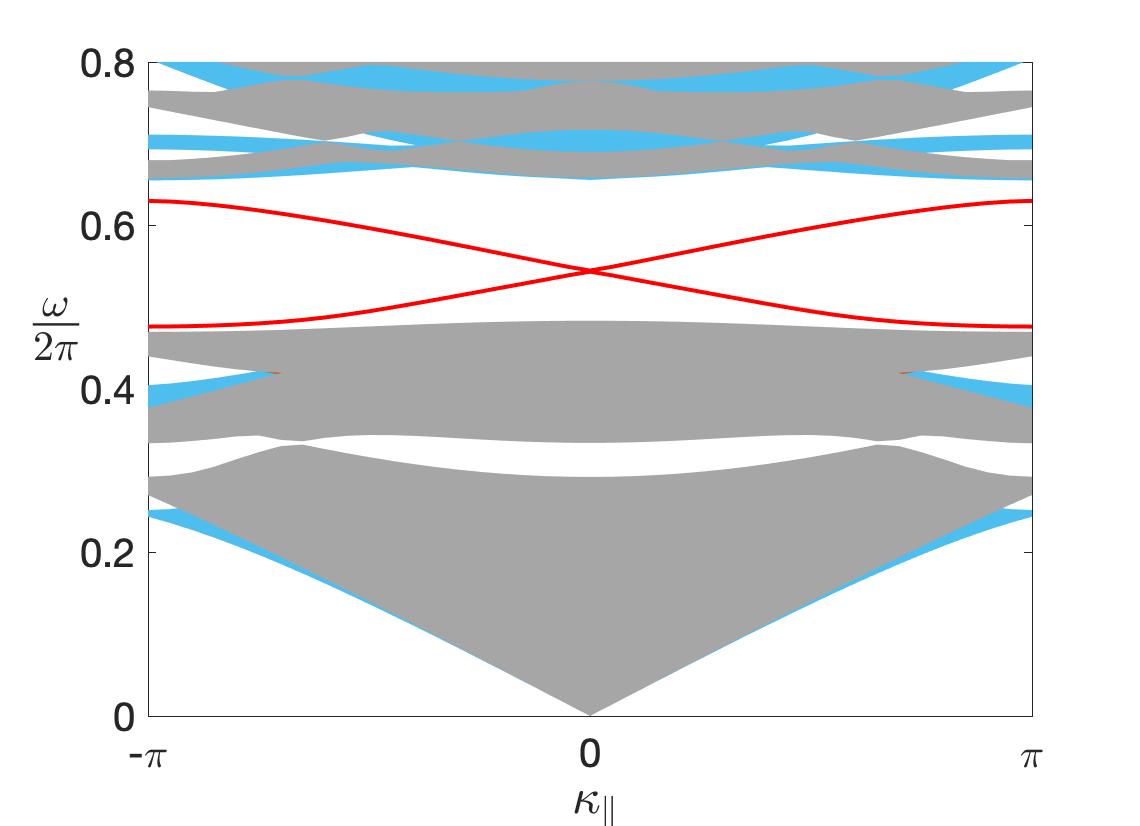}
\par\end{centering}
\caption{Top: The joint structure formed by two optimized periodic media that are shifted to the right with a distance of $0.025$; Bottom: The dispersion relation of the edge modes along the interface for the wavenumber $\kappa_\parallel\in[-\pi,\pi]$. Two dispersion curves cross with each other when $\kappa_\parallel=0$. See Section~\ref{s:ex3}. }
\label{Fig:tri_C6_edge_mode_dispersion2}
\end{figure}

Finally, we compute the bandwidth of edge modes for the joint structure by gluing the two optimized photonic crystals together along the interface direction $\be_2$. The dispersion relation of edge modes for $\kappa_\parallel\in[-\pi, \pi]$ is plotted in Figure \ref{Fig:tri_C6_edge_mode_dispersion1}. For each $\kappa_\parallel\in[-\pi, \pi]$, there exit two edge modes between the third and fourth band of the periodic media. It is clear that the bandwidth for edge modes is significantly enlarged after the optimization procedure. Indeed, the bandwidth of the edge modes for the joint structure can be further increased by shifting the periodic media before they are glued together. This is demonstrated in Figure \ref{Fig:tri_C6_edge_mode_dispersion2}, where the dispersion curves for the two edge modes cross with each other when $\kappa_\parallel=0$ if both media are shifted to the right with a distance of $0.025$. For completeness, we plot the edge modes when $\kappa_\parallel=0$ and $\pi$ in Row 3-6 of Figure \ref{Fig:tri_C6_edge_mode_dispersion1} for the joint structure shown in Row 1 of Figure \ref{Fig:tri_C6_edge_mode_dispersion1}. All the edge modes are localized along the interface of two media.

\section{Discussion} \label{s:Disc}
In this work, we present an efficient semi-definite programming method for
maximizing the shared spectral band-gap width of two topological photonic crystals. The topological invariants of the photonic crystals are preserved during the iterations.
The algorithm converges with a small number of iterations and the optimized media obtain significantly large band gap width, which allows for the transportation of wave energy in a large spectral gap via the topological edge modes between the two media. The computational framework can be extended to the TM polarized electromagnetic wave and the three-dimensional topological photonic crystals \cite{Ma_2016, oono2016section}. This will be reported in future work.

\section*{Appendix}
For an isolated band, a discrete Wilson loop along the $\bkappa_2$ direction for given $\bkappa_1$ is computed via the formula
\begin{equation}\label{eq:wilson_lopp1}
    W(\bkappa_1) = - \text{Im}\left(\ln \prod_{j=1}^{N-1} \langle \phi(\bkappa^j;\cdot), \phi(\bkappa^{j+1};\cdot) \rangle\right),
\end{equation}
where $\{\bkappa^j\}_{j=1}^N$ are a set of uniform grid points along the $\bkappa_2$ direction given by $\bkappa^j = \bkappa_1 + \frac{j}{N} \bkappa_2 $.
For a group of $m$ degenerate bands, let $\lambda (S)$ denote the eigenvalues of the matrix $S$, then
\begin{equation}\label{eq:wilson_lopp2}
    W(\kappa_x) = - \text{Im} \ln (\lambda (S)),
\end{equation}
where $S = \prod_{j=1}^{N-1} S_j$, and the $m\times m$ matrix $S_j$ is given by
\begin{equation*}
S_j =
\begin{bmatrix}
\langle \phi^{(1)} (\bkappa^{j};\cdot), \phi^{(1)} (\bkappa^{j+1};\cdot) \rangle  & \ldots & \langle \phi^{(1)} (\bkappa^{j};\cdot), \phi^{(m)} (\bkappa^{j+1};\cdot) \rangle  \\
\vdots & \ddots & \vdots \\
\langle \phi^{(m)} (\bkappa^{j};\cdot), \phi^{(1)} (\bkappa^{j+1};\cdot) \rangle  & \ldots & \langle \phi^{(m)} (\bkappa^{j};\cdot), \phi^{(m)} (\bkappa^{j+1};\cdot) \rangle
\end{bmatrix}.
\end{equation*}

\printbibliography

\end{document}